\documentclass[11pt]{article}

\usepackage{amsmath}
\usepackage{amsthm,amscd,amsfonts}
\usepackage{amssymb, upref, color}
\usepackage{amsmath,amssymb,amsthm,mathrsfs}
\usepackage{color}
\usepackage[dvips]{graphicx}
\usepackage{epsf,epsfig, verbatim} %subfigure,
\usepackage{latexsym,bm}
\usepackage{epsfig}
\usepackage{epstopdf}
\usepackage{subfig}
\usepackage[font=scriptsize]{caption}
\usepackage{graphicx}
\graphicspath{{figures/}}
\usepackage{booktabs}
\usepackage{float}
\usepackage{indentfirst}

\numberwithin{equation}{section}

\theoremstyle{plain}
\newtheorem{exam}{Example}[section]
\newtheorem{theorem}[exam]{Theorem}

\begin{document}
	\title{Turing instability in a diffusive predator-prey model with multiple Allee effect and herd behavior
\footnote{ This paper was supported by Zhejiang Provincial Natural Science Foundation of China (No.LZ23A010001), the NNSFC (No. 11671176 and 11931016).}}
	\author{Jianglong Xiao \,\,\,\,\, Yonghui Xia$\footnote{Corresponding author. Yonghui Xia, yhxia@zjnu.cn; xiadoc@163.com.}$
		\\
		{\small \textit{$^a$ College of Mathematics Science,  Zhejiang Normal University, 321004, Jinhua, China}}\\
		%{\small \textit{$^b$ Departamento de Matem\'aticas, Universidad de Chile, Santiago, Chile }}\\
		{\small Email: jianglongxiao@zjnu.edu.cn; yhxia@zjnu.cn; xiadoc@163.com.}
	}
	
	\maketitle

	\begin{abstract}
		 Diffusion-driven instability and bifurcation analysis are studied in a predator-prey model with herd behavior and  quadratic mortality by incorporating multiple Allee effect into prey species.
The existence and stability of the equilibria of the system are studied. The  sufficient and necessary conditions for Turing instability occurring are obtained. And the stability and direction of Hopf and steady state bifurcations are explored by using the normal form method.  Furthermore, some numerical simulations are presented to support our theoretical analysis. We found that too large diffusion rate of prey prevents Turing instability from emerging. The biomass conversion rate does affect the stability of the system and the occurrence of Turing instability. This indicates that the biomass conversion rate is essentially significant for the predator-prey system. Finally, we summarize our findings in the conclusion. \\
		{\bf keywords}: Turing instability; Hopf bifurcation; Steady state bifurcation;  Multiple Allee effect; Herd behavior;  Predator-prey;
	\end{abstract}
%{\bf MSC 2020} 34K23; 34D30; 37C60; 35B32; 35B35; 92D40
	
%%	{\bf  The study of the bifurcations and chaos of reaction-diffusion models from the ecosystems has been an interesting topic in recent years. In particular, Turing instability in diffusive models have been received great attentions from the biologist and mathematicians due to its great significance. In this paper, we investigated a  diffusive predator-prey model with multiple Allee effect, herd behavior. Moreover, we consider the quadratic mortality on predator species. We found that the dynamics of system near the positive equilibria is quite rich. We are more concerned with the spatial dynamics near the positive equilibrium. Turing instability is analyzed in the plane of the conversion rate of prey to predator and the diffusion rate of predator.  We found that too large diffusion rate of prey prevents Turing instability from emerging.  The biomass conversion rate does affect the stability of the equilibria and the occurrence of Turing instability. This indicates that the biomass conversion rate is greatly significant for the predator-prey system, and we can control the biological conversion rate to achieve the coexistence of the predator and prey. Finally, we summarize our findings in the conclusion.   }

	\section{Introduction}
	 \subsection{History}
	Modeling the interactions between the predator and prey in ecosystems using differential equations is one of the most popular methods in ecological research. And functional responses reflect the interactions between the predator and prey. Since the classical Lotka-Volterra model was built in Lotka \cite{Lotka} and Volterra \cite{Volterra}, substantial dynamic models with various functional responses have been put forward to study the relationship between the predator and prey. It's well known that there are some conventional functional responses, such as Holling  \uppercase\expandafter{\romannumeral1}-\uppercase\expandafter{\romannumeral4} types, Beddington-DeAngelis type, ratio-dependent  type and so on. In some predator-prey systems, prey exhibits herd behavior to defend themselves from predators and improve their survival ability. And there is indeed experimental evidence suggesting that not just prey, but predators exhibit schooling behavior, such as  Major \cite{P. F. Major},  Schmidt and Mech \cite{P. A. Schmidt}, Courchamp and Macdonald \cite{F. Courchamp}, Scheel and Packer \cite{D. Scheel}, Ajraldi et al. \cite{V. Ajraldi}. %For a predator-prey system in which prey exhibits herd behavior,
If prey exhibits herd behavior in the predator-prey system, then the square root  functional response proposed by Ajraldi et al.  \cite{V. Ajraldi} is more appropriate, and the model is formulated as follows:
	\begin{equation}\label{square root}
		\begin{cases}
			\frac{du}{dt}=ru(1-\frac{u}{K})-A\sqrt{u}v,\\
			\frac{dv}{dt}=B\sqrt{u}v-Dv,
		\end{cases}
	\end{equation}
    where the functions $u$ and $v$ denote the densities (at time $t$) of prey and predator, respectively; $r$ stands for the growth rate of prey; $K$  represents the maximal environmental carrying capacity; $A$ is regarded as the search efficiency of $v$ for  $u$; $B$ reflects the biomass conversion rate; $D$ is the death rate of the predator.
	Since then, many scholars have begun to study a predator-prey system with the square root functional response. Braza \cite{P. A. Braza} and Xu et al. \cite{C. Xu} further studied system \eqref{square root}. Tang and Song \cite{Tang and Song} studied the effect of cross-diffusion on system \eqref{square root}, their results show that  the cross-diffusion plays a considerable role in the pattern selection. Later, Tang et al. \cite{TangSongZhang} studied spatiotemporal dynamics %near Turing-Hopf bifurcation point
 of system \eqref{square root} with cross-diffusion, their results indicate that spatiotemporal dynamics %near 	the Turing-Hopf bifurcation point
is quite rich under proper conditions. From  the second formula of system \eqref{square root}, we see that $-Dv$ is regarded as linear mortality rate of the predator. In fact, a type of quadratic mortality was also considered, that is, we modified the term $-Dv$ to $-Dv^{2}$.
    Ghorai and  Poria  \cite{S. Ghorai} considered the effects of a diffusive predator-prey system with quadratic mortality rate and Holling \uppercase\expandafter{\romannumeral2} type. And Yuan et al. \cite{S. Yuan} investigated the  following  diffusive model:
	\begin{equation}\label{square root2}
		\begin{cases}
			\frac{\partial u}{\partial t}=d_{1}\Delta u+ru(1-\frac{u}{K})-A\sqrt{u}v,\\
			\frac{\partial v}{dt}=d_{2}\Delta v+B\sqrt{u}v-Dv^{2},
		\end{cases}
	\end{equation}
   where $d_{1}$ and $d_{2}$ are the diffusion coefficients of $u$ and $v$, respectively.
  Xu and Song  \cite{Xu and Song}  deliberated the system with herd behavior in order to distinctly explore the occurrence of Hopf bifurcation and Turing instability. And Singh and  Banerjee \cite{T. Singh} studied the system with herd behavior when  the mortality of the predator is linear or quadratic. For the system with herd behavior, Tang and Song \cite{TangSonghyperbolic} replaced the quadratic death rate with a more complicated hyperbolic death rate to study  bifurcation behaviors. Based on Tang and Song \cite{TangSonghyperbolic}, Tang et al. \cite{TangJiang} considered a delay effect on Hopf bifurcations. By incorporating the herd behavior into prey and schooling behavior into predator, Yang et al. \cite{YangYuanZhang} proposed  a novel functional response to study a predator-prey system with diffusion and delay. Moreover, Song and Tang \cite{SongTang} studied the system with herd behavior, prey-taxis and linear mortality. And Liu et al. \cite{LiuZhangMeng} investigated Turing-Hopf bifurcation of the system with herd behavior and prey-taxis. \\
	\indent In addition, Allee effect is one of the important phenomena affecting the population density of the predator or prey in ecology. And Allee effect has been observed in different organisms, such as plants, invertebrates and vertebrates in Berec et al.    \cite{L. Berec}. Generally speaking, there exist two kinds of Allee effects, one is the strong Allee effect, the other is the weak Allee effect. And Allee effect has been studied by lots of scholars in predator-prey systems. Wang et al. \cite{JingfengWang1} considered a  diffusive predator-prey system with the strong Allee effect. Later, Wang and Wei \cite{JingfengWang2} investigated a  diffusive and delayed predator-prey system with the strong Allee effect. Actually, a single population is even affected by two or more Allee effect simultaneously in Berec et al. \cite{L. Berec}. Pal  and Saha \cite{Pal Saha} qualitatively analyzed a predator-prey system with multiple Allee effect in prey. Singh et al. \cite{Singh Bhadauria} studied a modified Leslie-Gower predator-prey model with multiple Allee effect. And Tiwari and Raw \cite{Tiwari Raw} introduced the Crowley-Martin functional response into a diffusive Leslie-Gower predator-prey model. Feng and Kang \cite{Feng Kang} studied a modified Leslie-Gower model with multiple Allee effect.  And Wu et al. \cite{Wu Zhao Yuan} considered a diffusive predator-prey model with threshold harvesting and multiple Allee effect. Martinez and Aguirre \cite{Martinez Aguirre} studied a  Leslie-Gower predator-prey system with and Holling \uppercase\expandafter{\romannumeral1} type and multiple Allee effect. Based on Martinez and Aguirre \cite{Martinez Aguirre}, Li et al. \cite{LiYangWeiWang} investigated Turing-Hopf bifurcation. Tiwari et al. \cite{Tiwari Raw Mishra} studied the multiple Allee effect on a prey-predator model with schooling behaviour, their spatiotemporal patterns show that the strong Allee effect drives the populations to develop high-density schools. Not just the predator-prey system, but there are other types of diffusive systems, such as Lengyel-Epstein system studied by Yi et al. \cite{Fengqi Yi1}. And Yi et al. \cite {Fengqi Yi2} also studied  spatiotemporal patterns and bifurcations in a  diffusive predator-prey system. Huang and Wu \cite{Huang Wu} considered a stage-structured SLIRM epidemic model with latent period. And Wu et al. \cite{Wu Yang Zou} investigated the spatial-temporal dynamics in a competitive system with nonlocal dispersals.
	\subsection{Motivation and model formulation}
	The above works indicate that the study of multiple Allee effect is of great practical significance to the predator-prey system. However, there is no study taking into account the factors of the multiple Allee effect, herd behavior, and quadratic mortality together until now. Therefore, based on system \eqref{square root2}, we propose the following model under the homogeneous Neumann boundary condition:
	\begin{equation}\label{motivation system1}
		\begin{cases}
			\frac{\partial u}{\partial t}=D_{1}\Delta u+\frac{ru}{u+N}(1-\frac{u}{K})(u-M)-A\sqrt{u}v, \ \  \ x\in \Omega,t \textgreater 0,\\
			\frac{\partial v}{\partial t}=D_{2}\Delta v+B\sqrt{u}v-Dv^2, \ \ \ \ \ \ \ \ \ \ \ \ \ \ \ \ \ \ \ \ \ \ \ \  x\in\Omega,t\textgreater0,\\
			 u_\upsilon= v_\upsilon=0, \ \ \ \ \ \ \ \ \ \ \ \ \ \ \ \ \ \ \ \ \ \ \ \ \ \ \ \ \ \ \ \ \ \ \ \ \ \ \ \ \ \ \ \ \ t \textgreater 0,\\
			u(0,x)=u_{0}(x) \geq 0,v(0,x)=v_{0}(x) \geq 0, \ \ \ \ \ \ \ \ \ x \in \Omega,
		\end{cases}
	\end{equation}
	where $M$ denotes Allee threshold with $-K\textless M \textless K$, $N$ is viewed as the positive auxiliary parameter, $\Delta$ is the Laplacian operator.  $\Omega=(0,\pi)$ and $\upsilon$ account for a bounded domain with smooth boundary and the outward unit normal vector to $\Omega$, respectively. $D_{1}$ and $D_{2}$ are the diffusion coefficients of $u$ and $v$, respectively. And $\frac{ru}{u+N}$ stands for another Allee effect due to external factors affecting prey birth rate, $N$ satisfies $N \textgreater -M$. Allee effect is weak if $-K\textless M \textless 0$, Allee effect is strong if $0\textless M \textless K$.\\
	 By  transformation
	\begin{equation*}
		\stackrel{-}{u}=\frac{u}{K}, \stackrel{-}{v}=\frac{A}{r \sqrt{K}}v, \stackrel{-}{t}=rt, a=\frac{B \sqrt{K}}{r}, b=\frac{D \sqrt{K}}{A}, m=\frac{M}{K}, n=\frac{N}{K}, d_{1}=\frac{D_{1}}{r}, d_{2}=\frac{D_{2}}{r},
	\end{equation*}
	and after dropping the bars, system \eqref{motivation system1} is rescaled to
	\begin{equation}\label{diffusion system}
		\begin{cases}
			\frac{\partial u}{\partial t}=d_{1}\Delta u+\frac{u}{u+n}(1-u)(u-m)-\sqrt{u}v, \ \ \ \  \ x\in \Omega,t \textgreater 0,\\
			\frac{\partial v}{\partial t}=d_{2}\Delta v+\theta v(\sqrt{u}-cv),  \ \ \ \ \ \ \ \ \ \ \ \ \ \ \ \ \ \ \ \ \ \ \  x\in\Omega,t\textgreater0,\\
		\end{cases}
	\end{equation}
where $\theta=a, c=\frac{b}{\theta}$, the boundary and initial conditions are the same as in system \eqref{motivation system1}.\\
	 This paper is devoted to studying the existence and stability of equilibria, Turing instability and bifurcation  behaviors of system \eqref{diffusion system} with or without diffusion.
	\subsection{Organization of paper}
	The existence and local stability of equilibria are discussed  and some bifurcations of system \eqref{diffusion system} without diffusion are given in Section 2. Turing instability, Hopf bifurcation and steady state bifurcation of system \eqref{diffusion system} are explored in  Section 3.  Some numerical simulations and calculations are carried out to support our theoretical results in Section 4. Finally, we end this paper with a summary of our findings.  We found that too large diffusion rate $d_{1}$ of prey prevents Turing instability from emerging.

	\section{Analysis of the local system (the system without diffusion)}
	For reaction-diffusion system \eqref{diffusion system}, the local system has the following form:
	\begin{equation}\label{local system}
		\begin{cases}
			\frac{\partial u}{\partial t}=\frac{u}{u+n}(1-u)(u-m)-\sqrt{u}v,\\
			\frac{\partial v}{\partial t}=\theta v(\sqrt{u}-cv).
		\end{cases}
	\end{equation}
	
	\subsection{Equilibria}
	In this section, we investigate the existence of  all equilibria. Obviously, system \eqref{local system} has three boundary equilibria $E_{0}(0,0)$, $E_{1}(1,0)$ and $E_{2}(m,0)$. Since the system cannot be linearized at $E_{0}(0,0)$, we don't study the stability of $E_{0}(0,0)$ here. And all  equilibria satisfy
	\begin{equation*}
		\begin{cases}
			\frac{u}{u+n}(1-u)(u-m)-\sqrt{u}v=0,\\
			\theta v(\sqrt{u}-cv)=0.
		\end{cases}
	\end{equation*}
	By simple calculation, all possible  positive equilibria $(u,v)$ satisfy $cu^{2}+(1-(m+1)c)u+mc+n=0$ and $v=\frac{\sqrt{u}}{c}$. \\
	Let $\Delta =\left(1-(m+1)c\right)^{2}-4c(mc+n)$ and let $m_{1}=-mc, m_{2}=\frac{\left(1-(m+1)c\right)^{2}}{4c}-mc$. Then we have:
	\begin{theorem}\label{equilibria theorem}
		System \eqref{diffusion system} always has three boundary equilibria $E_{0}(0,0)$, $E_{1}(1,0)$ and $E_{2}(m,0)$. Moreover, for interior equilibrium, we have:\\
		(a) If $ n \textless m_{1}$, then there exists a unique positive equilibrium $E_{30}(u_{30},v_{30})$, where \\
		$u_{30}=\frac{(m+1)c-1+\sqrt{\Delta}}{2c}, v_{30}=\frac{(m+1)c-1+\sqrt{\Delta}}{2c^{2}}$.\\
		(b) If $n = m_{1}$ and $(m+1)c \textgreater 1$, then there exists a positive equilibrium $E_{31}(u_{31},v_{31})$, where \\
		$u_{31}=\frac{(m+1)c-1}{c}, v_{31}=\frac{(m+1)c-1}{c^{2}}$.\\
		(c) If $m_{1} \textless n \textless m_{2}$ and $(m+1)c \textgreater 1$, then there exist two positive equilibria $E_{32}(u_{32},v_{32}), E_{30}(u_{30},v_{30})$,\\
		where 	$u_{32}=\frac{(m+1)c-1-\sqrt{\Delta}}{2c},  v_{32}=\frac{(m+1)c-1-\sqrt{\Delta}}{2c^{2}}$, \\
		(d) If $n=m_{2}$ and $(m+1)c \textgreater 1$, then there exists a unique positive equilibrium $E_{33}(u_{33},v_{33})$, where $u_{33}=\frac{(m+1)c-1}{2c},  v_{33}=\frac{(m+1)c-1}{2c^{2}}$.\\
		(e) If  $n \textgreater m_{2}$, then there is no  positive equilibrium.
	\end{theorem}
	\noindent $\bm{Proof}.$
	Let $h(u)=cu^{2}+(1-(m+1)c)u+mc+n=0$, then the symmetry is $u=\frac{(m+1)c-1}{2c}$.\\
	If $n \textless m_{1}$, then $\Delta \textgreater 0$ and $h(0) \textless 0$, which implies system \eqref{local system} has a positive equilibrium $E_{30}(u_{30},v_{30})$.\\
	If $n = m_{1}$ and $(m+1)c \textgreater 1$, then $h(0)=0$ and the symmetry $u=\frac{(m+1)c-1}{2c} \textgreater 0$, which means that system \eqref{local system} has a positive equilibrium $E_{31}(u_{31},v_{31})$.\\
	If $m_{1} \textless n \textless m_{2}$ and $(m+1)c \textgreater 1$, then $\Delta \textgreater 0$ and $h(0) \textgreater 0$, which means that system \eqref{local system} has two positive equilibria $E_{32}(u_{32},v_{32}), E_{30}(u_{30},v_{30})$.\\
	If $n=m_{2}$ and $(m+1)c \textgreater 1$, then $\Delta =0$ and the symmetry $u=\frac{(m+1)c-1}{2c} \textgreater 0$, which suggests that system \eqref{local system} has a positive equilibrium $E_{33}(u_{33},v_{33})$.\\		
	If  $n \textgreater m_{2}$, then $\Delta \textless0$, which suggests that system \eqref{local system} has no positive equilibrium.
	
	\subsection{Local stability}
	In this subsection, we study the local stability of equilibria. In view of the length of the paper, we only list the results here and put detailed proof of theorem \ref{E1 theorem}-\ref{E33 theorem} in  the appendix 1.
	\begin{theorem}\label{E1 theorem}
		The boundary equilibrium $E_{1}(1,0)$ is  a saddle for $-1 \textless m\textless 1$.
	\end{theorem}
	
	\begin{theorem}\label{E2 theorem}
		The boundary equilibrium $E_{2}(m,0)$ is an unstable  node for $0 \textless m\textless 1$.
	\end{theorem}
	
	\begin{theorem}\label{E30 theorem}
		For the  positive equilibrium $E_{30}(u_{30},v_{30})$, we have:\\
		(i) If  $\theta \textgreater \theta_{30}$, then $E_{30}(u_{30},v_{30})$ is asymptotically stable;\\
		(ii) If $\theta \textless \theta_{30}$, then $E_{31}(u_{30},v_{30})$ is unstable.
	\end{theorem}

	\begin{theorem}\label{E31 theorem}
		For the  positive equilibrium $E_{31}(u_{31},v_{31})$, we have:\\
		(i) If  $\theta \textgreater \theta_{31}$, then $E_{31}(u_{31},v_{31})$ is asymptotically stable;\\
		(ii) If $\theta \textless \theta_{31}$, then $E_{31}(u_{31},v_{31})$ is unstable.
	\end{theorem}

    \begin{theorem}\label{E32 theorem}
       The  positive equilibrium $E_{32}(u_{32},v_{32})$ is a saddle.
    \end{theorem}
	
	\begin{theorem}\label{E33 theorem}
		For the  positive equilibrium $E_{33}(u_{33},v_{33})$, we have:\\
		(a) If  $\theta \neq \theta_{33}$, then $E_{33}$ is a saddle node; \\
		(b) If  $\theta = \theta_{33}$,  then \\
		\indent  (i) $E_{33}$ is a cusp of codimension at least 3 if  $2\alpha_{2}+\alpha_{3}=0$;\\
		\indent   (ii) $E_{33}$ is a cusp of codimension 2 if  $2\alpha_{2}+\alpha_{3}\neq0$.\\
	\end{theorem}

	\subsection{Bifurcation analysis of the local system}
	In this subsection, we study the transcritical bifurcation, the saddle-node  bifurcation and Hopf bifurcation of system \eqref{local system}. In view of the length of the paper, we put detailed proof of theorem \ref{Transcritical theorem1}-\ref{Saddle-node theorem} in  the appendix 2.
	\subsubsection{Transcritical bifurcation}
	From  theorem \ref{E1 theorem}, we observe that $m=1$ turns $E_{1}(1,0)$ to a nonhyperbolic equilibrium. Thus,  system  \eqref{local system} may exhibit the  transcritical bifurcation at $E_{1}$.
	\begin{theorem}\label{Transcritical theorem1}
		If the intensity of Allee effect arrives maximum (at $m=1$ ), then  system \eqref{local system} exhibts
		the transcritical bifurcation at $E_{1}$.  And the  bifurcation parameter  is $m=m_{TC}=1$.
	\end{theorem}

	\begin{theorem}\label{Transcritical theorem2}
		If the intensity of  Allee effect arrives maximum (at $m=1$ ), then  system \eqref{local system} exhibts
		the transcritical bifurcation at $E_{2}$.  And the  bifurcation parameter  is $m=m_{TC}=1$.
	\end{theorem}
	
	\subsubsection{Saddle-node  bifurcation}
	From theorem \ref{equilibria theorem}, we observe that the number of the interior equilibrium of system \eqref{local system} alters  when parameter $n$ passes from one side of $n=n_{SN}$ to the other side. Thus,  system  \eqref{local system} may exhibit the   saddle-node  bifurcation at $E_{33}(u_{33},v_{33})$.
	\begin{theorem}\label{Saddle-node theorem}
		System \eqref{local system} exhibts the  saddle-node  bifurcation at $E_{33}$ if $n=m_{2}$, $(m+1)c \textgreater 1$ and $\theta\neq \theta_{33}$.  And the  bifurcation parameter  is $n=n_{SN}=\frac{((m+1)-1)^{2}}{4c}-mc$.
	\end{theorem}

	\subsubsection{Hopf bifurcation}
	Here we regard $\theta$ as the bifurcation parameter to study Hopf bifurcation of system \eqref{local system}. According to theorem \ref{E30 theorem}- \ref{E31 theorem} and $\frac{d\text{Tr}(E_{3i})}{d\theta}=-\sqrt{u_{3i}} \textless 0(i=0,1)$, we easily get
	\begin{theorem}
		If $\theta = \theta_{30}$, then system \eqref{local system} exhibts Hopf bifurcation at $E_{30}(u_{30},v_{30})$.
	\end{theorem}
		
	\begin{theorem}
	     If $\theta = \theta_{31}$, then system \eqref{local system} exhibts Hopf bifurcation at $E_{31}(u_{31},v_{31})$.
	\end{theorem}

	\section{Analysis of the diffusive system}
	\subsection{Turing instability and bifurcation analysis}
	In this subsection, we analyze Turing instability of positive equilibrium  $E_{31}(u_{31},v_{31})$ and the existence of Hopf  bifurcation and steady state bifurcation. In the following, we restrict ourselves to
	the  case: $\Omega =(0, \pi)$. For convenience, we denote $\frac{(1-(m+1)c)u_{31}}{c(u_{31}-mc)}+\frac{1}{2c}$ by $\delta_{1}$ and denote $\sqrt{u_{31}}$ by $\delta_{2}$.
	Linearizing system \eqref{diffusion system} at $E_{31}(u_{31},v_{31})$ , we get
	\begin{equation}
		\left(
		\begin{array}{c}
			\frac{\partial u }{\partial t}\\\frac{\partial v }{\partial t}
		\end{array}
		\right)
		:= L
		\left(
		\begin{array}{c}
			u \\ v
		\end{array}
		\right)
		=D\Delta
		\left(
		\begin{array}{c}
			u \\  v
		\end{array}
		\right)
		+J
		\left(
		\begin{array}{c}
			u \\ v
		\end{array}
		\right),
	\end{equation}
	where
	\begin{equation}
		D=
		\left(
		\begin{array}{cc}
			d_{1}\Delta & 0\\
			0     & d_{2}\Delta \\
		\end{array}
		\right),
		J=
		\left(
		\begin{array}{cc}
			\delta _{1} & -\delta _{2} \\
			\frac{\theta}{2c} & -\theta \delta _{2}
		\end{array}
		\right).
	\end{equation}
	Define the real-valued Sobolev space
	\begin{equation*}
		X=\{(u,v) \in W^{2,2}(0,\pi)|u_{x}(0,t)=u_{x}(\pi,t)=0,v_{x}(0,t)=v_{x}(\pi,t)=0\}
	\end{equation*}
	with the inner product for $u=(u_{1},u_{2})^{T},v=(v_{1},v_{2})^{T} \in X$
	\begin{equation}
		[u,v]=\sum_{i=1}^{2} \int_{0}^{\pi} u_{i} v_{i} dx.
	\end{equation}
    It is known to all that  the eigenvalue problem:
    \begin{equation*}
    -\ddot{\phi}=\mu \phi, \ x \in \Omega; \ \dot{\phi}(0)=\dot{\phi}(\pi)=0,
    \end{equation*}
     has eigenvalues $\mu_{k}=k^2(k\in\mathbb{N})$ with eigenfunctions $\varepsilon_{k}(x)=\frac{\cos kx}{||\cos kx||_{2,2}}$.
	On the other hand, the eigenvalues of $L$ are given by the eigenvalues of $L_{k}$, and the characteristic equation of $L_{k}:=J-\text{diag}\{d_{1}k^{2},d_{2}k^{2}\}$  is
	\begin{equation}\label{quadratic equations}
		\lambda ^{2}-T_{k}\lambda+D_{k}=0,\ k\in \mathbb {N},
	\end{equation}
	where
	\begin{flalign}
		&T_{k}=\delta _{1}-\theta \delta _{2}-(d_{1}+d_{2})k^{2},\\
		&D_{k}=d_{1}d_{2}k^{4}+(d_{1} \theta \delta _{2}-d_{2} \delta _{1})k^{2}+\frac{\theta \delta _{2}}{2c}(1-2c\delta _{1}).
	\end{flalign}
	In what follows, we restrict $\delta_{1}\textgreater0$.
	According to $T_{k}=0$, we obtain
	\begin{equation}\label{Hopf curve}
		d_{2}=d_{2}^{H}(k,\theta) \triangleq -\frac{\delta _{2}}{k^{2}}\theta +\frac{\delta _{1}-d_{1}k^{2}}{k^{2}}.
	\end{equation}
	According to $D_{k}=0$, we obtain
	\begin{equation}\label{Turing curve}
		d_{2}=d_{2}^{T}(k,\theta) \triangleq \frac{\delta _{2}(d_{1}k^{2}+\frac{1}{2c}-\delta _{1})}{(\delta _{1}-d_{1}k^{2})k^{2}} \theta.
	\end{equation}
	Substituting \eqref{Hopf curve} into \eqref{Turing curve} yields
	\begin{equation}
		\theta =\theta _{k}^{*} \triangleq \frac{2c(\delta _{1}-d_{1}k^{2})^{2}}{\delta _{2}}.
	\end{equation}
	Solving \eqref{Hopf curve} with respect to $\theta$ yields
	\begin{equation}\label{True Hopf curve}
		\theta=\theta ^{H}(k,d_{2}) \triangleq -\frac{d_{2}}{\delta _{2}}k^{2}+\frac{\delta _{1}-d_{1}k^{2}}{\delta _{2}}, \ 0 \leq k \leq k^{*},
	\end{equation}
	where $k^{*}=\mathrm{max} \{k\in \mathbb{N}| \delta_{1}-d_{1}k^{2}\textgreater 0 \}$. Thus, for fixed $k\in [0,k^{*}]$, equation \eqref{quadratic equations} has pairs of pure imaginary roots $\pm i\sqrt{D_{k}}$ if $d_{2}=d_{2}^{H}(k,\theta)$ and $\theta \textgreater \theta _{k}^{*}$.\\
	Solving \eqref{Turing curve} with respect to $\theta$ yields
	\begin{equation}\label{True Turing curve}
		\theta=\theta ^{T}(k,d_{2}) \triangleq \frac{(\delta _{1}-d_{1}k^{2})k^{2}}{d_{1}\delta _{2}k^{2}+\delta _{2}(\frac{1}{2c}-\delta _{1})}d_{2}.
	\end{equation}
	The curve determined by $D_{k}=0$ in the $\theta - d_{2}$ plane is called the Tuing bifurcation curve and denoted by $l_{k}$. \\
	Choosing $\theta$ as a parameter and  denoting the root of \eqref{quadratic equations} by
	$\lambda (\theta)$. Taking the derivative of both sides of \eqref{quadratic equations} with respect to $\theta$ results
	\begin{equation*}
		2\lambda \frac{d \lambda (\theta)}{d\theta}-\lambda \frac{dT_{k}(\theta)}{d\theta}-T_{k}\frac{d \lambda (\theta)}{d\theta}+\frac{dD_{k}(\theta)}{d\theta}=0.
	\end{equation*}
	Then we get that the following transversality conditions hold,
	\begin{flalign}
		&\frac{d Re \lambda (\theta)}{d \theta} \Big|_{\theta=\theta ^{H}(k,d_{2})}=-\frac{\delta_{2}}{2} \textless 0,\\
		&\frac{d Re \lambda (\theta)}{d \theta} \Big|_{\theta=\theta ^{T}(k,d_{2})}=\frac{d_{1}\delta _{2}k^{2}+\delta _{2}(\frac{1}{2c}-\delta _{1})}{T_{k}}\textless 0,
	\end{flalign}
	which implies that system \eqref{diffusion system} exhibits steady state bifurcation at $E_{31}(u_{31},v_{31})$ on $l_{k}$.
	Define the following curve in the $\theta-d_{2}$ plane,
	\begin{equation}\label{H_{k}}
		H_{k}:d_{2}=d_{2}^{H}(k,\theta), \theta \textgreater \theta_{k}^{*}, \ k=0,1,\cdots, k^{*},
	\end{equation}
	which represents  Hopf bifurcation curve of system \eqref{diffusion system} at positive equilibrium $E_{31}(u_{31},v_{31})$. In particular, $H_{0}$ represents Hopf bifurcation curve  of the local system.
	In addition, for the local system, $E_{31}(u_{31},v_{31})$ is asymptotically stable if  $\theta \textgreater \frac{\delta_{1}}{\delta_{2}}$ and $E_{*}(u_{*},v_{*})$ is unstable if  $\theta \textless \frac{\delta_{1}}{\delta_{2}}$. \\
	\indent Based on the above analysis, we have:
	\begin{theorem}\label{Hopf and steady state}
		Suppose  $n=m_{1}$ and $(m+1)c \textgreater 1$ hold. Let $H_{k}$ and $d_{2}=d_{2}^{T}(k,\theta)$  in the $\theta-d_{2}$ plane be defined by \eqref{H_{k}} and \eqref{Turing curve}, respectively. Then we get:\\
		(a) System \eqref{diffusion system} exhibits the steady state bifurcation at
		$E_{*}(u_{*},v_{*})$ on the line $d_{2}=d_{2}^{T}(k,\theta)$.\\
		(b) System \eqref{diffusion system} exhibits Hopf bifurcation at $E_{31}(u_{31},v_{31})$ when $\theta=\theta ^{H}(k,d_{2})$. Spatially homogeneous periodic solutions arise on $H_{0}$, and spatially inhomogeneous periodic solutions arise on $H_{k}$, $k=  1,\cdots, k^{*}$.		
	\end{theorem}
	Next we study  Turing instability of the  equilibrium $E_{31}(u_{31},v_{31})$ of system \eqref{diffusion system}.  The first quadrant of the $\theta-d_{2}$ plane is split into two parts by the line $H_{0}$, and 0-mode Hopf bifurcation curve  is always located above $k$-mode ($1 \leq k \leq k_{*}$) Hopf bifurcation curve. This indicates the corresponding solution of $k$-mode Hopf bifurcation is always unstable. That is, we only study Turing instability above $H_{0}$, which satifies $\theta  \textgreater \frac{\delta_{1}}{\delta_{2}}$. Thus, we have:\\
	\begin{theorem}
		Assume that $n=m_{1}$, $(m+1)c \textgreater 1$ and $\frac{\delta_{1}}{\delta_{2}} \textless \theta \textless \theta ^{T}(k,d_{2}) \ ( \text{for some} \  k \in \mathbb{N}^{+})$ hold, then Turing instability occurs, where $\theta ^{T}(k,d_{2})$ is defined by \eqref{True Turing curve}.
	\end{theorem}
 Additionally, we denote the slope  of $l_{k}$ as $\eta_{k}=\frac{(\delta _{1}-d_{1}k^{2})k^{2}}{d_{1}\delta _{2}k^{2}+\delta _{2}(\frac{1}{2c}-\delta _{1})}$. Then $\eta_{k}\textgreater 0$ if $1\leq k \leq k^{*}$ and $\eta_{k}\textless 0$ if $k \textgreater k^{*}$. Apparently, $l_{k}$ does not intersect $H_{0}$. This indicates that when $k \textgreater k^{*}$, Turing-Hopf bifurcation does not occur in system \eqref{diffusion system}. Hence, we restrict our discussion to the case $1\leq k \leq k^{*}$. To consider the monotonicity of $\eta_{k}$  with respect to $k$, let
	\begin{equation}
		\eta(x)=\frac{(\delta _{1}-d_{1}x)x}{d_{1}\delta _{2}x+\delta _{2}(\frac{1}{2c}-\delta _{1})}, \ x \in [1,k^{*}].
	\end{equation}
	By taking the derivative, we have that  $\eta(x)$ maximizes at the point $x^{*}$, where
	\begin{equation}
		x^{*}=\frac{2c\delta _{1}-1+\sqrt{1-2c\delta _{1}}}{2cd_{1}}.
	\end{equation}
	Let
	\begin{equation}
		k_{m}=
		\begin{cases}
			[\sqrt{x^{*}}], \ \ \ \ \ \ \eta \left([\sqrt{x^{*}}]+1\right) \leq \eta \left([\sqrt{x^{*}}]\right),\\
			[\sqrt{x^{*}}]+1, \ \eta \left([\sqrt{x^{*}}]+1\right) \textgreater \eta \left([\sqrt{x^{*}}]\right),
		\end{cases}
	\end{equation}
	where $[.]$ represents the integer part function.
	Thus, there exists a positive integer $k_{m}$ such that when $k=k_{m}$, $\eta(k)$ takes the maximum value. That is, $H_{0}$ intersects  $l_{k_{m}}$ at $(d_{2}^{m},\theta _{m})$. And a $(k_{m},0)$-mode Turing-Hopf bifurcation occurs at the intersection $(d_{2},\theta)=(d_{2}^{m},\theta _{m})$. Moreover, the line $l_{k_{m}}$ splits the stability region into two parts, one is a Turing unstable  region $T$ and the other is a still stable region $S$. \\
	\indent In order to verify the above theoretical analysis, we take the following parameters:\\
	(H1): $d_{1}=0.01,d_{2}=0.02,c=\frac{100}{41},m=-0.5,n=\frac{50}{41}$;\\
	(H2): $d_{1}=0.1,c=\frac{100}{41},m=-0.5,n=\frac{50}{41}$.\\
	When system \eqref{diffusion system} takes the parameters of (H1), we have $\delta_{1}=0.1988, \delta_{2}=0.3,\theta=0.6627$ for $T_{0}=0$ and $\theta^{T}=0.7777$ for wave number $k=1$. And we know that $E_{31}(u_{31},v_{31})$ is locally asymptotically stable by theorem \ref{E31 theorem}. Thus, we restrict the parameter $\theta \in(0.6627,0.7777)$ in order for Turing instability to occur under conditions (H1), see Fig. 1.\\
	On the other hand, when system \eqref{diffusion system} takes the parameters of (H2), we have $k^{*}=1$ and $\eta(x)=\frac{-0.1x^{2}+0.1988x}{0.03x+0.0018}$, see $(c)$ of Fig.2. Then we get $k_{m}=1$ and the following curves:\\
	$H_{0}: \theta=0.6627, H_{1}: \theta=-3.3333d_{2}+0.3294, L_{1}: \theta=3.1019d_{2}$,\\
	where $H_{0}$ and $H_{1}$ stand for Hopf bifurcation curve, $l_{1}$ stands for Turing bifurcation curve. Then we plot the  Turing-Hopf bifurcation curves in the $\theta-d_{2}$ plane. And Turing bifurcation curve $l_{1}$ intersects with Hopf bifurcation curve $H_{0}$ at the point $(d^{m}_2,\theta_{m}) =(0.2136,0.6627)$, see $(b)$ of Fig. 2.
	\begin{figure}[H]
		\centering
		\begin{minipage}[c]{0.5\textwidth}
			\centering
			\includegraphics[height=5cm,width=8cm]{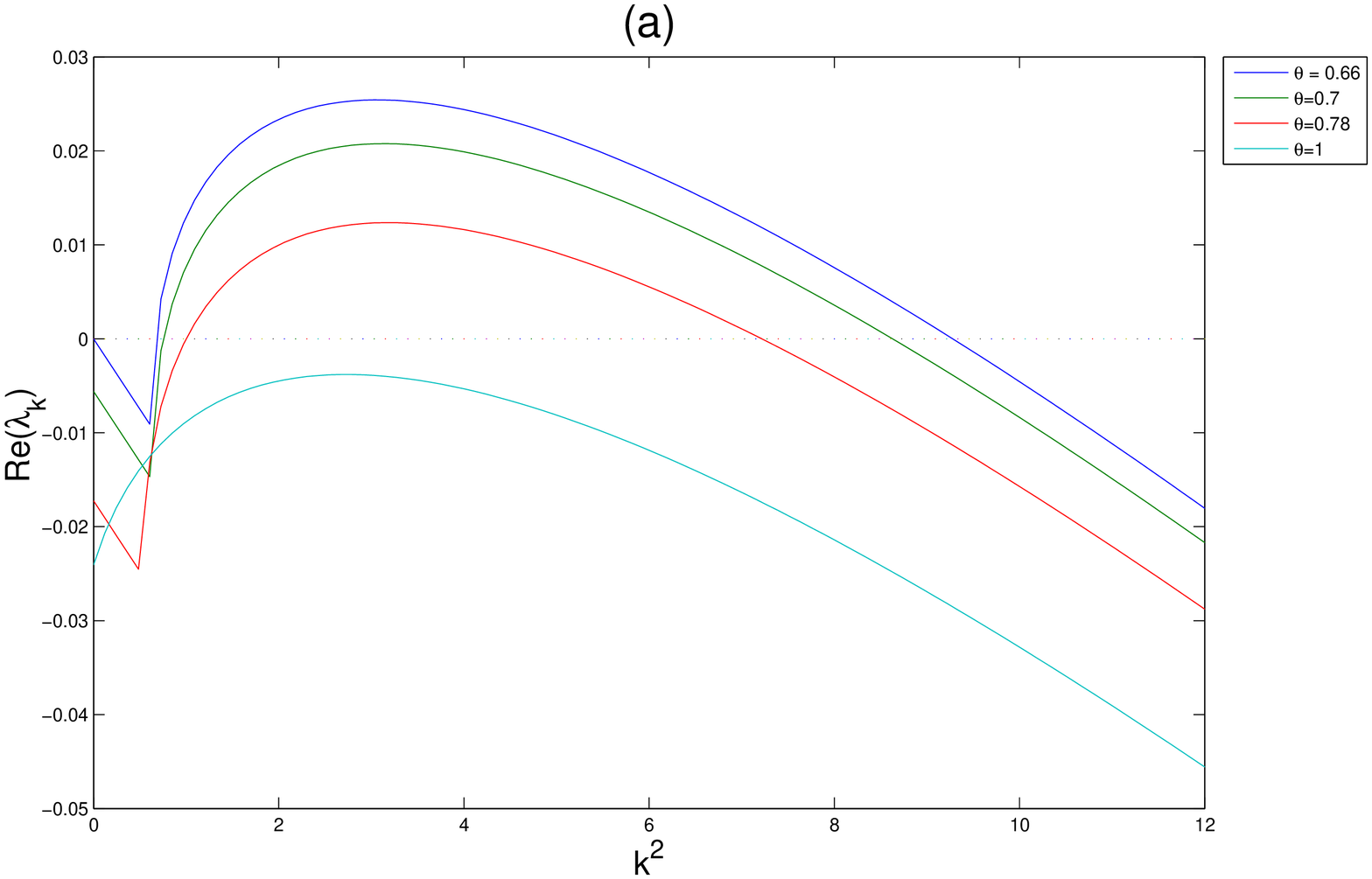}
		\end{minipage}%
		\begin{minipage}[c]{0.5\textwidth}
			\centering
			\includegraphics[height=5cm,width=8cm]{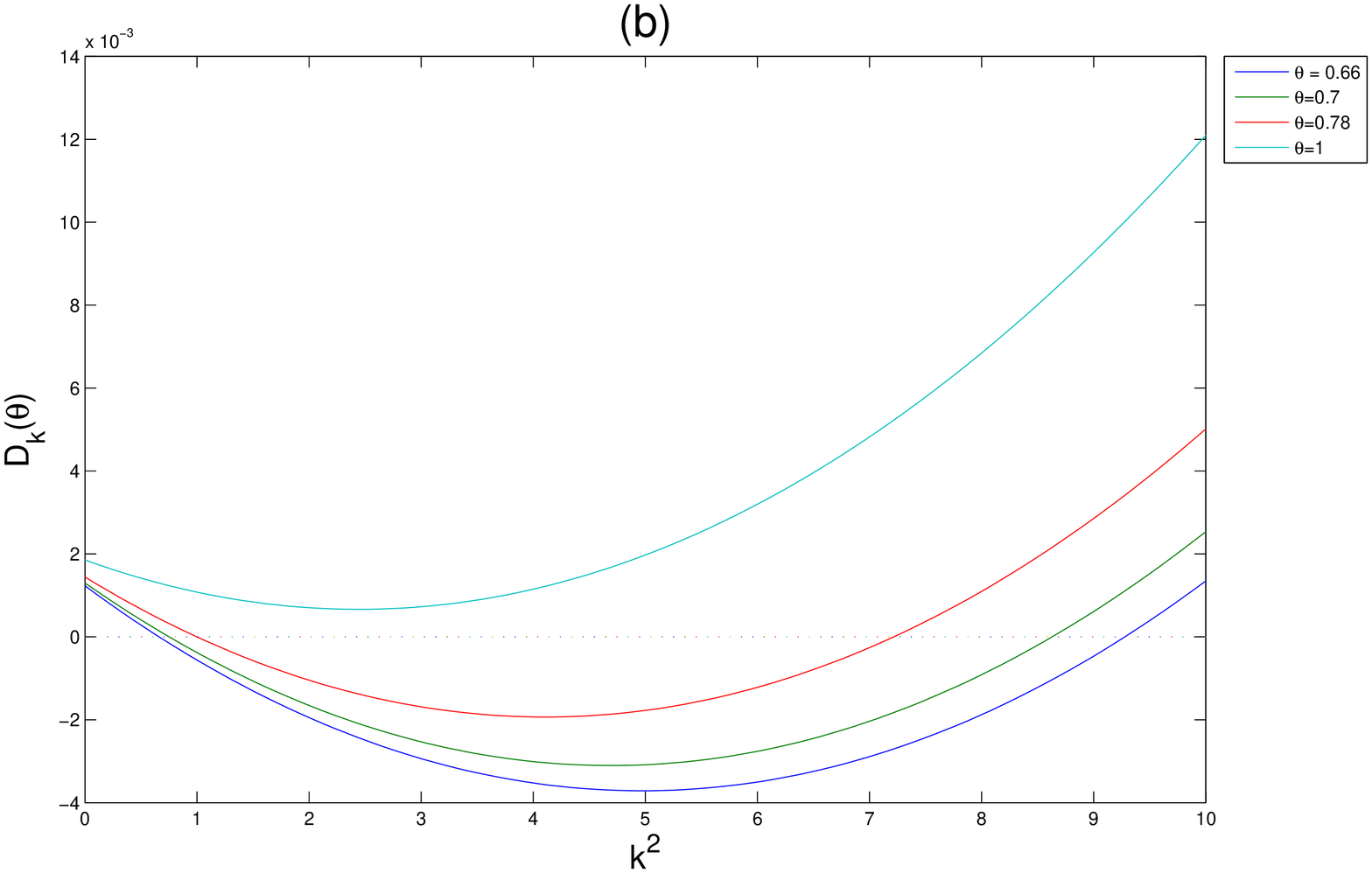}
		\end{minipage}
		\caption{(a) Dispersion relation diagram with respect to $k^2$. (b) Graph of $D_{k}(\theta)$ with respect to  $k^2$ for different $\theta.$}
	\end{figure}
	
	\begin{figure}[H]
		\centering
		\begin{minipage}[c]{0.5\textwidth}
			\centering
			\includegraphics[height=5cm,width=8cm]{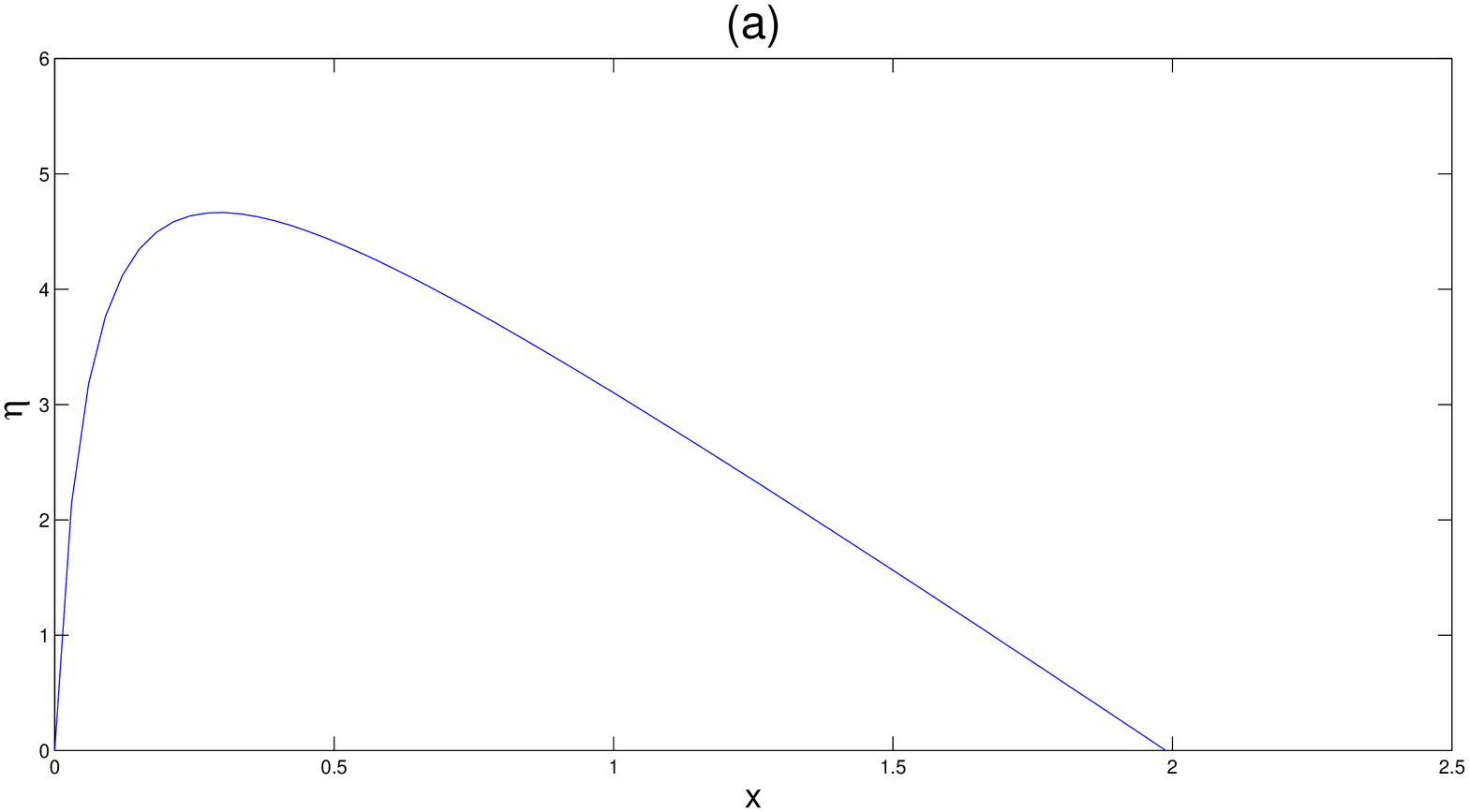}
		\end{minipage}%
		\begin{minipage}[c]{0.5\textwidth}
			\centering
			\includegraphics[height=5cm,width=8cm]{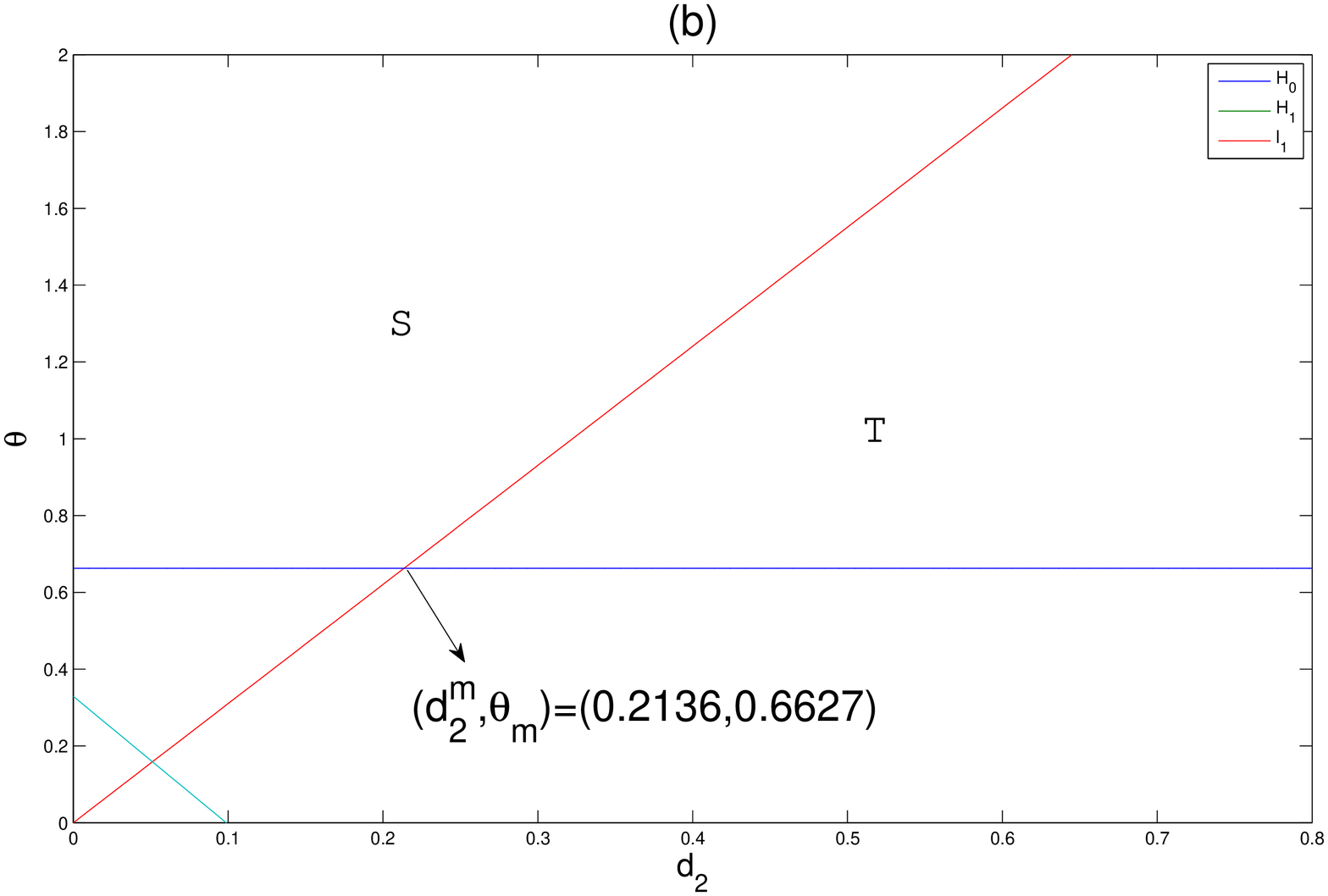}
		\end{minipage}
		\caption{(a) Graph of $\eta(x)=\frac{-0.1x^{2}+0.1988x}{0.03x+0.0018}$. (b) Turing-Hopf bifurcation diagram in the $\theta-d_{2}$ plane.}
	\end{figure}

	\subsection{Normal form for Hopf bifurcation and steady state bifurcation}
   In this subsection, we use  the center manifold and the normal form theory
   to judge the the direction and the stability of Hopf bifurcation and steady state bifurcation. Choosing $\theta$ as a parameter and denoting the critical value by $\theta_{*}$. And $\theta _{*}$ = $\theta ^{H}(k,d_{2})$ and $\theta _{*}$ = $\theta ^{H}(k,d_{2})$ are defined by \eqref{True Hopf curve} for Hopf bifurcation and \eqref{True Turing curve} for steady state bifurcation, respectively. We introduce a perturbation parameter	$\epsilon$ by setting  $\epsilon=\theta -\theta_{*}$, then  $\epsilon=0$ is the bifurcation value. And we rewrite  $E_{31}(u_{31},v_{31})$ as $E_{31}(u_{31}(\epsilon),v_{31}(\epsilon))$. Let $\hat {u}(.,t)={u}(.,t)-u_{31}(\epsilon), \hat {v}(.,t)={v}(.,t)-v_{31}(\epsilon)$, and after dropping the bars, we get
	\begin{equation}\label{partial system}
		\frac{dU(t)}{dt}=d\Delta U(t)+L_{0}(U(t))+G(U(t),\epsilon),
	\end{equation}
	where
	\begin{equation}
		U(t)=
		\left(
		\begin{array}{c}
			u(t) \\ v(t)
		\end{array}
		\right),
		d\Delta=
		\left(
		\begin{array}{cc}
			d_{1}\frac{{\partial}^{2}}{\partial x^{2}} & 0\\
			0     & d_{2}\frac{{\partial}^{2}}{\partial x^{2}}  \\
		\end{array}
		\right),
		L_{0}(U(t))=
		\left(
		\begin{array}{c}
			\delta_{1}u-\delta_{2}v \\ \frac{\theta \delta_{1}}{2c}u-\theta \delta_{2}v
		\end{array}
		\right),
	\end{equation}
	$$G(U(t),\epsilon)=\sum _{i+j+s\geq 2}\frac{1}{i!j!s!}f_{ijs}u^{i}v^{j}\epsilon^{s},\ 	f_{ijs}=\left(f^{(1)}_{ijs},\ f^{(2)}_{ijs}\right)^{T},$$
	with $f^{(k)}_{ijs}=\frac{\partial^{i+j+s}}{\partial u^{i} \partial v^{j} \partial \epsilon^{s}}\hat{f}^{(k)}(0,0,0), \ k=1,2,$ and
	\begin{flalign*}
		&\hat{f}^{(1)}(u,v,\epsilon)=\frac{(u+u_{31}(\epsilon)) (1-(u+u_{31}(\epsilon)) (u+u_{31}(\epsilon)-m)}{u+u_{31}(\epsilon)+n}-\sqrt{u+u_{31}(\epsilon)}(v+v_{31}(\epsilon)),\\
		&\hat{f}^{(2)}(u,v,\epsilon)=(\theta_{*}+\epsilon)(v+v_{31}(\epsilon))\left(\sqrt{u+u_{31}(\epsilon)}-c(v+v_{31}(\epsilon))\right).
	\end{flalign*}
	The linearized system of \eqref{partial system} at the origin follows that
	\begin{equation}\label{linear partial differential equation}
		\frac{d U(t)}{dt}=\mathcal {L}(U(t)),
	\end{equation}
	where $\mathcal {L}(U(t))=d\Delta U(t)+L_{0}(U(t))$.
	Let $\mathcal{B}_{k}=span \{[\psi(.),\zeta ^{j}_{k}] \zeta ^{j}_{k}| \psi \in X, j=,1,2 \}$, then we have $L_{0}(\mathcal{B}_{k})\subset span\{ \zeta ^{1}_{k},\zeta ^{2}_{k} \} \ ,k\in \mathbb{N}^{+}$. Suppose
	$y(t)\in \mathbb{R}^{2}$ and
	\begin{equation}
		y^{T}(t)
		\left(
		\begin{array}{c}
			\zeta ^{1}_{k} \\ \zeta ^{2}_{k}
		\end{array}
		\right)
		\in \mathcal{B}_{k}.
	\end{equation}
	Then, on $\mathcal{B}_{k}$,
	 \eqref{linear partial differential equation} is equivalent to the following equation  on $\mathbb{R}^{2}$
	\begin{equation}\label{ linear ODE}
		\dot {y}(t)=
		\left(
		\begin{array}{cc}
			-d_{1}k^{2} & 0 \\
			0     &    	-d_{2}k^{2}
		\end{array}
		\right)
		y(t)+L_{0}(y(t)),
	\end{equation}
	where $y(t)\in \mathbb{R}^{2}$.
	Obviously, \eqref{ linear ODE} and \eqref{partial system} have the same characteristic equation \eqref{quadratic equations}.\\
	Let
	\begin{equation}
		\mathcal{H}_{k}=
		\left(
		\begin{array}{cc}
			\delta_{1}-d_{1}k^{2} & -\delta_{2} \\
			\frac{\theta_{*}}{2c}   &    	-\theta_{*}\delta_{2} -d_{2}k^{2}
		\end{array}
		\right)
	\end{equation}
	be the  characteristic matrix of \eqref{ linear ODE}. Denote the finite set of all eigenvalues of $\mathcal{H}_{k}$ with zero real parts by $\Lambda_{k}$. Then $\Lambda_{k}$ is also the  finite set of all eigenvalues of \eqref{linear partial differential equation} with zero real parts.\\
	\indent Next, we compute the normal forms for Hopf bifurcation and steady
	state bifurcation as per the method in Song et al. \cite{SongZhangPeng} and  Song \cite{SongZou}, from which we can judge the stability and direction of Hopf bifurcation and pitchfork bifurcation.
	
	\subsubsection{Stability  and  direction of Hopf bifurcation}
	For Hopf bifurcation, we have $\Lambda_{k}=\{ \text{i} \omega_{k}, -\text{i} \omega_{k}\}, R_{k}=\text{diag}\{\text{i} \omega_{k}, -\text{i} \omega_{k}\}, p=2$.\\
	\begin{theorem} \label{Hopf bifurcation theorem}
		Assume that there is a $s\in \mathbb{N}$ such that the characteristic equation \eqref{quadratic equations} with $\theta=\theta_{*}=\theta^{H}(s)$  has  pairs of pure imaginary roots $\pm \text{i} \omega_{s}(=\pm \text{i}\sqrt{D_{s}})$  and the remaining roots of  \eqref{quadratic equations} have nonzero real parts. Consquently,\\
		(i) if $v_{s2}\textless 0$, then the Hopf bifurcation at the
		critical value $\theta=\theta^{H}(s)$ is supercritical, and the bifurcating periodic solution is asymptotically stable;\\
		(ii) if $v_{s2}\textgreater 0$, then the Hopf bifurcation at the
		critical value $\theta=\theta^{H}(s)$ is subcritical, and the bifurcating periodic solution is unstable.
	\end{theorem}
	\noindent$\bm{Proof}.$
	Let
	\begin{equation*}
		p_{s}=
		\left(
		\begin{array}{c}
			\frac{2c(d_{2}s^{2}+\delta_{2}\theta_{*}+\text{i} \omega_{s})}{\theta_{*}}\\  1
		\end{array}
		\right)
		\triangleq
		\left(
		\begin{array}{c}
			p_{s1}\\ p_{s2}
		\end{array}
		\right),
		q_{s}=
		\left(
		\begin{array}{c}
			\frac{\theta_{*}}{4ic \omega_{s}}\\  \frac{d_{1}s^{2}-\delta_{1}+\text{i}\omega_{s}}{2\text{i}\omega_{s}}
		\end{array}
		\right)
		\triangleq
		\left(
		\begin{array}{c}
			q_{s1}\\ q_{s2}
		\end{array}
		\right).
	\end{equation*}
	Predictably, $\mathcal{H}_{s}p_{s}=\text{i}\omega_{s}p_{s},\mathcal{H}^{T}_{s}q_{s}=\text{i}\omega_{s}q_{s}$ and $<q^{T}_{s},p_{s}>=1$. \\
	On the basis of $\Phi_{s}=(p_{s},\bar{p}_{s}), \Psi_{s}=col(q^{T}_{s},\bar{q}^{T}_{s})$, then $<\Phi_{s},\Psi_{s}>=\Phi_{s}\Psi_{s}=I_{2}$, where $I_{2}$ is the identity matrix of order 2. Now  we  decompose $(u,v)^{T}$ as follows
	\begin{equation}
		\left(
		\begin{array}{c}
			u \\v
		\end{array}
		\right)
		=(z_{1}p_{s}+z_{2}\bar p_{s}) \varepsilon_{s}(x)+\omega,
	\end{equation}
	where $ z_{1}, z_{2}\in \mathbb{R}, \omega=({\omega_{1},\omega_{2}})^{T}$.
	By Song et al. \cite{SongZhangPeng}, Song and Zou \cite{SongZou}, we obtain that for Hopf bifurcation, the normal form truncated to the third terms takes form as
	\begin{equation}\label{norm form }
		\dot{z}=R_{s}z+
		\left(
		\begin{array}{c}
			R_{s1}z_{1}\varepsilon \\  \overline R_{s2}z_{2}\varepsilon
		\end{array}
		\right)
		+
		\left(
		\begin{array}{c}
			R_{s2}z^{2}_{1}z_{2}\varepsilon \\  \overline R_{s2}z_{1}z^{2}_{2}\varepsilon
		\end{array}
		\right)
		+O(|z| |\varepsilon|^{2}+|z|^{4}),
	\end{equation}
	where
	\begin{equation}
		R_{s1}=(f^{(2)}_{101} p_{s1}+f^{(2)}_{011} p_{s2})q_{s2}
	\end{equation}
	and
	\begin{equation}
		R_{s2}=
		\begin{cases}
			\frac{1}{2\pi}b_{021}+\frac{1}{4\pi}c_{021}+\frac{1}{2\sqrt{\pi}}E_{(0,0)}, \ s=0,\\
			\frac{3}{4\pi}b_{s21}+\frac{1}{2\sqrt{\pi}}E_{(s,0)}+\frac{1}{2\sqrt{2\pi}}E_{(s,2s)}, \ s\neq0, \\
		\end{cases}
	\end{equation}
	where \\
	$b_{s21}={q}^{T}_{s}\left(f_{300}p_{s1}|p_{s1}|^{2}+f_{030}p_{s2}|p_{s2}|^{2}+f_{210}(p^{2}_{s1} \bar{p}_{s2}+2p_{s2}|p_{s1}|^{2})+f_{120}(p^{2}_{s2} \bar{p}_{s2}+2p_{s1}|p_{s2}|^{2})\right)$,\\
	$c_{s21}=\frac{\text{i}}{w_{s}}\left(({q}^{T}_{s}A_{s20})({q}^{T}_{s}A_{k11})-|{q}^{T}_{s}A_{s11}|^{2}-\frac{2}{3}|{q}^{T}_{s}A_{s02}|^{2} \right)$  \\
	with \\
	\begin{flalign*}
		A_{s20}&=\overline {A}_{s02}=f_{200}p^{2}_{s1}+2f_{110}p_{s1}p_{s2}+f_{020}p^{2}_{s2},\\ A_{s11}&=f_{200}p^{2}_{s1}+4f_{110}\mathbf{Re} \{p_{s1}\bar{p}_{s2}\}+2f_{020}|p^{2}_{s2}|^{2},
	\end{flalign*}
	and \\
	\begin{flalign*}
		E_{(s,j)}&={q}^{T}_{s}\left((f_{200}p_{21}+f_{110}p_{s2})h^{(1)}_{sj11}+(f_{110}p_{s1}+f_{020}p_{s2})h^{(2)}_{sj11}+(f_{200}\bar p_{21}+f_{110}\bar p_{s2})h^{(1)}_{sj20}\right. \\
		&+\left.(f_{110}\bar p_{21}+f_{020}\bar p_{s2})h^{(2)}_{sj20}\right)
	\end{flalign*}
	with
	\begin{flalign*}
		h_{0020}=&\frac{1}{\sqrt{\pi}}(2\text{i}w_{0}I_{2}-\mathcal{H}_{0})^{-1}(A_{020}-q^{T}_{0}A_{020}p_{0}-\bar q^{T}_{0}A_{020}\bar p_{0}),\\
		h_{0011}=&-\frac{1}{\sqrt{\pi}}\mathcal{H}_{0}^{-1}(A_{011}-\bar q^{T}_{0}A_{011}\bar p_{0}-\bar q^{T}_{0}A_{011}\bar p_{0}),\\
		h_{sj20}=&\sigma_{sj}(2\text{i}w_{s}I_{2}-\mathcal{H}_{j})^{-1}A_{s20},\ s\neq 0, \ j=0, 2s, \\
		h_{sj11}=&-\sigma_{sj}\mathcal{H}_{j}^{-1}A_{s11},\ s\neq 0, \ j=0, 2s, \\
	\end{flalign*}
	and
	\begin{equation*}
		\sigma_{sj}=\int_{0}^{\pi}\varepsilon^{2}_{s}(x)\varepsilon_{j}(x)dx=
		\begin{cases}
			\frac{1}{\sqrt{\pi}}, \ \ j=0,\\
			\frac{1}{\sqrt{2\pi}}, \ \ j=2s\neq0,\\
			0,\ \ \rm otherwise.
		\end{cases}
	\end{equation*}
	By substituting variables $z_{1}=v_{1}-\text{i}v_{2}, z_{2}=v_{1}+\text{i}v_{2}$, the norm form \eqref{norm form }  can be transformed to the real coordinates. And then changing to cylindrical coordinates through $v_{1}=\rho cos \chi, v_{2}=\rho sin\chi$, namely
	\begin{equation}
		\begin{cases}
			\dot \rho=&v_{s1}\epsilon \rho+v_{s2}\rho^{3}+O(\epsilon \rho^{2}+|(\epsilon,\rho)|^{4}),\\
			\dot\chi=&-\omega_{s}+O(|(\epsilon,\rho)|),
		\end{cases}
	\end{equation}
	where $v_{s1}=\mathbf{Re}\{R_{s1}\}, v_{s2}=\mathbf{Re}\{R_{s2}\}$.
	By consulting Wiggins \cite{Wiggins}, we know that if $v_{s1}v_{s2}\neq0 $, then the direction of the bifurcation and the stability of the nontrivial periodic orbits are  determined by the sign of $v_{s2}$. Therefore,  $v_{s2}\textless 0$ suggests that a supercritical and stable Hopf bifurcation at the threshold value $\theta=\theta^{H}(s)$ occurs.  $v_{s2}\textgreater 0$ suggests that a subcritical and unstable  Hopf bifurcation at the threshold value $\theta=\theta^{H}(s)$ occurs.

	\subsubsection{Stability  and  direction of pitchfork bifurcation}
	For pitchfork bifurcation, we have $\Lambda_{k}=\{ 0\}, R_{k}=0, p=1$.\\
	\begin{theorem}\label{pitchfork bifurcation theorem}
		Assume that  there is a positive
		integer $s\in \mathbb{N}^{+}$  such that the characteristic equation \eqref{quadratic equations} with $\theta=\theta_{*}=\theta^{T}(s)$  has  a simple zero root $\lambda = 0 $ and the remaining roots of \eqref{quadratic equations}  have nonzero real parts. Consquently,\\
		(i) if $Q_{s30}\textless 0$, then system \eqref{diffusion system} exhibits a supercritical pitchfork bifurcation around $E_{31}(u_{31},v_{31})$ at the critical value $\theta=\theta^{T}(k)$;\\
		(ii) if $Q_{s30}\textgreater 0$, then system \eqref{diffusion system} exhibits a subcritical pitchfork bifurcation around $E_{31}(u_{31},v_{31})$ at the critical value $\theta=\theta^{T}(s)$.
	\end{theorem}
	\noindent$\bm{Proof}.$
	Let
	\begin{equation*}
		\widetilde{p} _{s}=
		\left(
		\begin{array}{c}
			1  \\  \frac{\delta_{1}-d_{1}s^{2}}{\delta_{2}}
		\end{array}
		\right)
		\triangleq
		\left(
		\begin{array}{c}
			\widetilde{p}_{s1}\\ \widetilde{p}_{s2}
		\end{array}
		\right),
		\widetilde{q}_{s}=
		\left(
		\begin{array}{c}
			-\frac{d_{2}s^{2}+\delta_{2}\theta_{*}}{T_{s}}\\  \frac{\delta_{2}}{T_{s}}
		\end{array}
		\right)
		\triangleq
		\left(
		\begin{array}{c}
			\widetilde{q}_{s1}\\ \widetilde{q}_{s2}
		\end{array}
		\right).
	\end{equation*}
	Predictably, $\mathcal{H}_{s}\widetilde{p}_{s}=\text{i}\omega_{s}\widetilde{p}_{s},\mathcal{H}^{T}_{k}\widetilde{q}_{s}=\text{i}\omega_{s}\widetilde{q}_{s}$ and $<\widetilde{q}^{T}_{s},\widetilde{p}_{s}>=1$. \\
	On the basis of $\Phi_{s}=\widetilde{p}_{s}, \Psi_{s}=\widetilde{q}^{T}_{s}$, then $<\Phi_{s},\Psi_{s}>=\Phi_{s}\Psi_{s}=I_{2}$. Now  we  decompose $(u,v)^{T}$ as follows:
	\begin{equation}
		\left(
		\begin{array}{c}
			u \\v
		\end{array}
		\right)
		=\widetilde{p}_{s}z \varepsilon_{s}(x)+\omega,
	\end{equation}
	where $ z\in \mathbb{R}, \omega=({\omega_{1},\omega_{2}})$.\\
	By Song et al. \cite{SongZhangPeng}, Song and Zou \cite{SongZou}, we obtain that for the steady state bifurcation, the normal form truncated to third terms  takes forms as
	\begin{equation}
		\dot{z} = Q_{s11}\epsilon z +Q_{s30}z^{3},
	\end{equation}
	where $Q_{s11}=(f^{(2)}_{101} \widetilde{p}_{s1}+f^{(2)}_{011} \widetilde{p}_{s2})\widetilde{q}_{s2}$ and $Q_{s30}=\frac{1}{4\pi}\gamma_{s}+\frac{1}{2\sqrt{\pi}}\gamma_{(s,0)}+\frac{1}{2\sqrt{2\pi}}\gamma_{(s,2s)}$ \\
	with
	$\gamma_{s}=\widetilde{q}^{T}_{s}\left(f_{300}\widetilde{p}_{s1}^{3}+f_{030}\widetilde{p}_{s2}^{3}+3f_{210}\widetilde{p}^{2}_{s1} \widetilde{p}_{s2}+3f_{120}\widetilde{p}_{s1}\widetilde{p}^{2}_{s2} \right)$,\\
	$\gamma_{(s,j)}=\widetilde{q}^{T}_{s}\left((f_{200}\widetilde{p}_{s1}+f_{110}\widetilde{p}_{s2})h^{(1)}_{sj}+(f_{110}\widetilde{p}_{s1} +f_{020}\widetilde{p}_{s2})h^{(2)}_{sj} \right), \ j=0, \ 2s$, \\
	with $h_{sj}=\varepsilon_{sj}\mathcal{H}^{-1}_{j}A_{s20}$.\\
	By consulting Wiggins \cite{Wiggins}, we know that if $Q_{s30}Q_{s11} \neq 0 $, then the direction of  steady state bifurcation and the stability  are  determined by the sign of $Q_{s30}$. Therefore, $Q_{s30}\textless 0$ suggests that a supercritical pitchfork bifurcation around $E_{31}(u_{31},v_{31})$ at the threshold value $\theta=\theta^{T}(s)$ arises.
	$Q_{s30}\textgreater 0$ suggests that a subcritical pitchfork bifurcation around $E_{31}(u_{31},v_{31})$ at the threshold value $\theta=\theta^{T}(s)$ arises.

		\section{ Numerical simulations}
\begin{figure}[H]
		\centering
		\begin{minipage}[c]{0.5\textwidth}
			\centering
			\includegraphics[height=5cm,width=7cm]{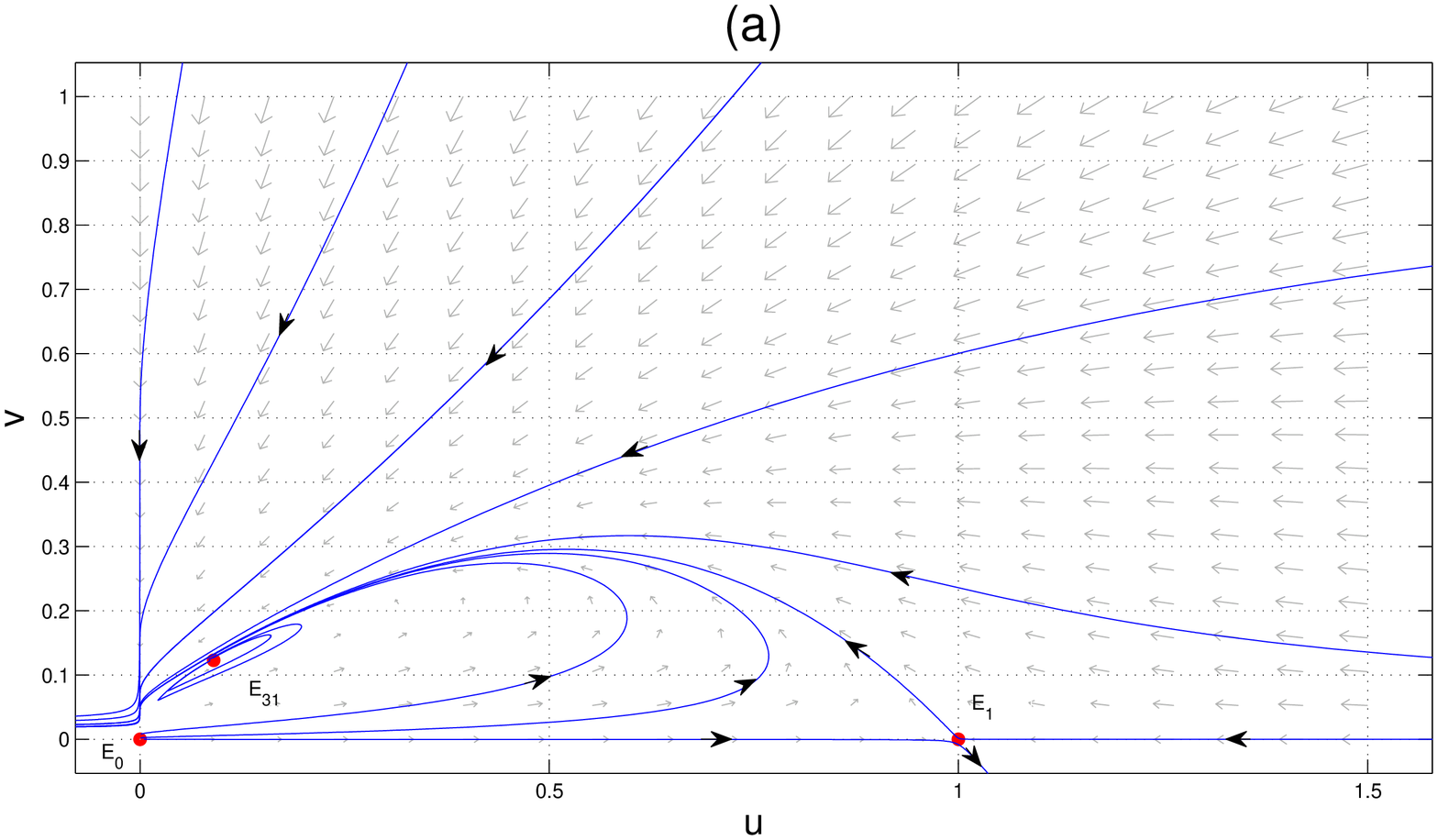}
		\end{minipage}%
		\begin{minipage}[c]{0.5\textwidth}
			\centering
			\includegraphics[height=5cm,width=7cm]{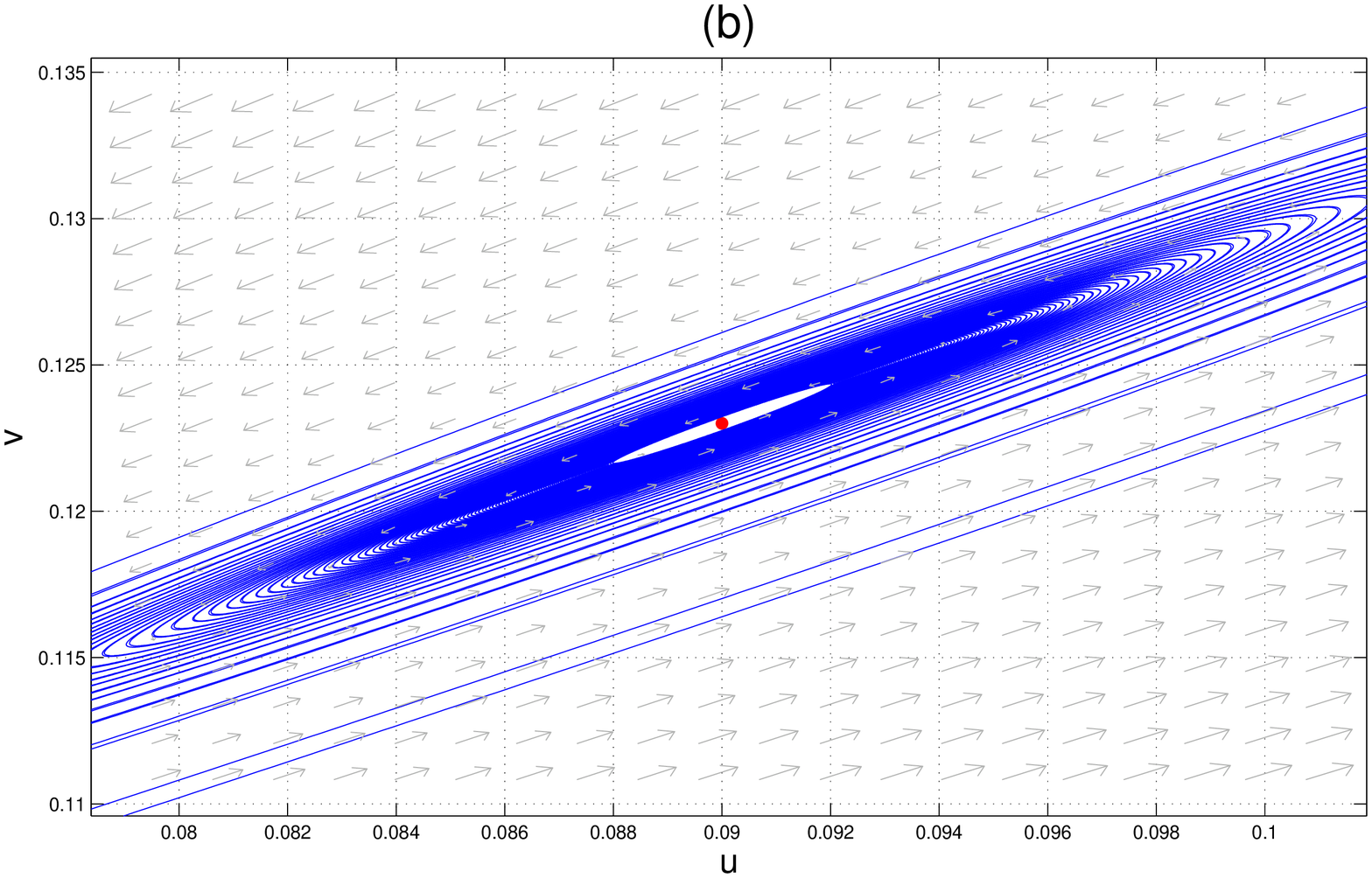}
		\end{minipage}
	\caption{(a) $E_{31}(0.09,0.123)$ is asymptotically stable if $\theta=0.68 \textgreater \theta_{0}=0.6627$. (b) There is a stable limit cycle arising from  $E_{31}$ if $\theta=0.662\textless \theta_{0}=0.6627$.}\label{localE31stable}
	\end{figure}

\begin{figure}[H]
		\centering
		\begin{minipage}[c]{0.5\textwidth}
			\centering
			\includegraphics[height=5cm,width=7cm]{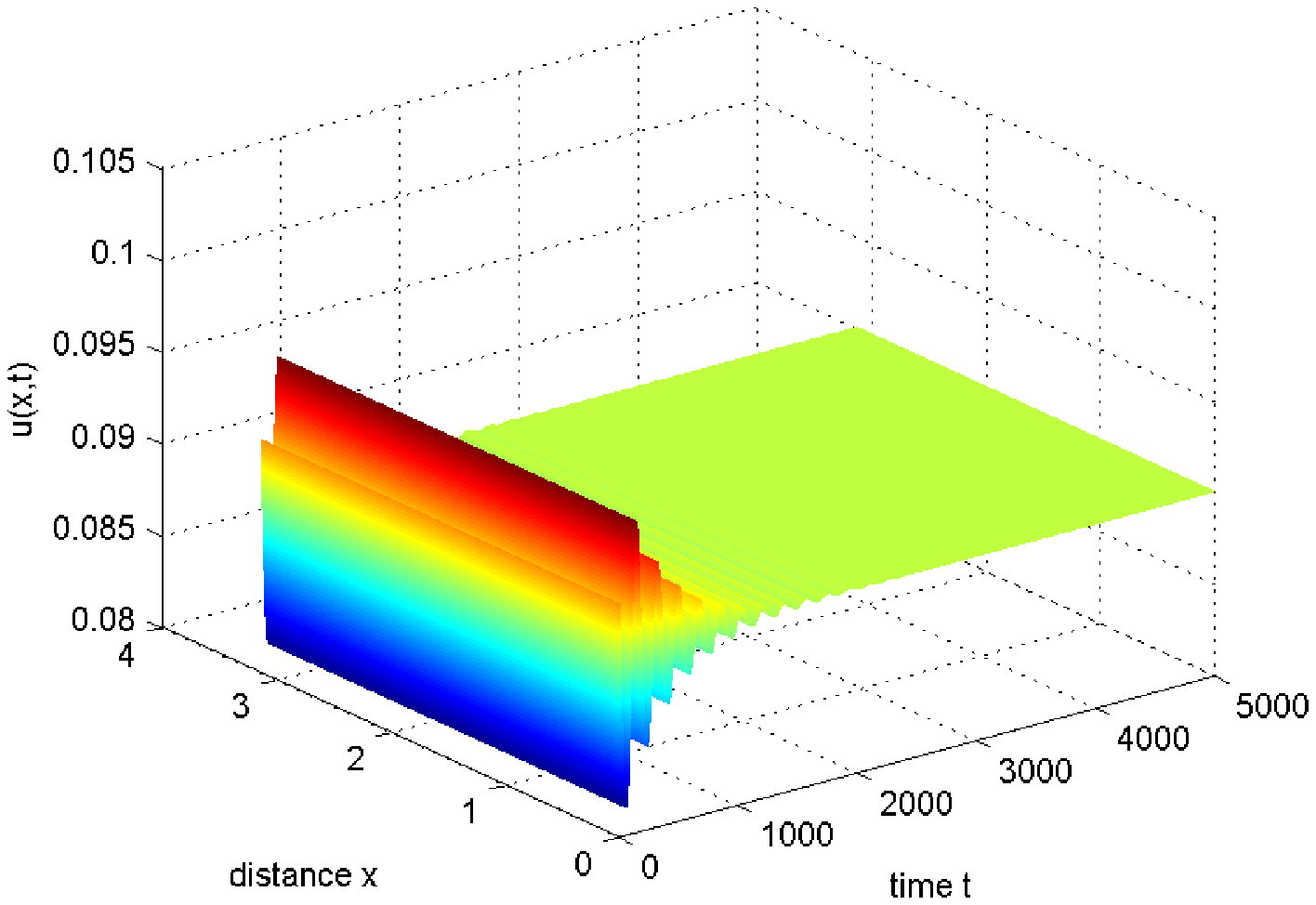}
		\end{minipage}%
		\begin{minipage}[c]{0.5\textwidth}
			\centering
			\includegraphics[height=5cm,width=7cm]{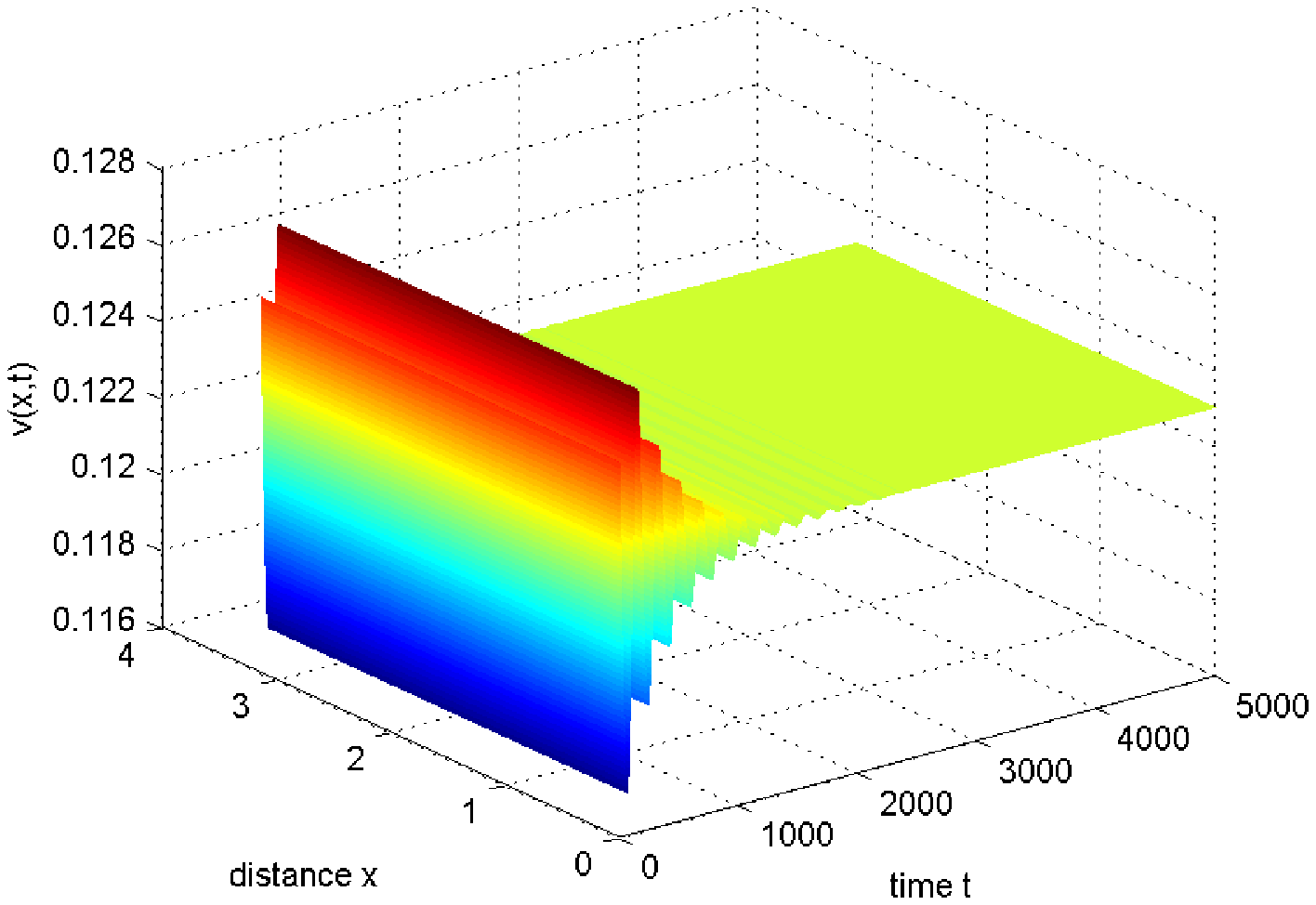}
		\end{minipage}
		\caption{For fixed parameter values  $d_{1}=0.1, d_{2}=0.4,m=-0.5,n=\frac{50}{41},c=\frac{100}{41},\theta=0.68\textgreater \theta_{0}=0.6627$, the positive equilibrium  $E_{31}(0.09,0.123)$ is asymptotically stable. The initial values are $(u{(x,0)},v{(x,0)})=(0.093,0.126)$.}\label{diffusionE31stable}
	\end{figure}

	In this section, we provide some numerical illustrations to support aforementioned analysis. Continue with the boundary and initial conditions in system \eqref{motivation system1} and consider the following system:
	\begin{equation}\label{Numerical diffusion system}
		\begin{cases}
			\frac{\partial u}{\partial t}=0.1\Delta u+\frac{u}{u+\frac{50}{41}}(1-u)(u+0.5)-\sqrt{u}v, \ \ \  \ x\in \Omega,t \textgreater 0,\\
			\frac{\partial v}{\partial t}=d_{2}\Delta v+\theta v(\sqrt{u}-\frac{100}{41}v),  \ \ \ \ \ \ \ \ \ \ \ \ \ \ \ \ \ \ \ \ \ \ \  x\in\Omega,t\textgreater0.\\
		\end{cases}
	\end{equation}
	In system \eqref{Numerical diffusion system}, we take parameters in (H2): $d_{1}=0.1,c=\frac{100}{41},m=-0.5,n=\frac{50}{41}$. Then we get that system \eqref{Numerical diffusion system} has a unique positive equilibrium $E_{31}(0.09,0.123)$ by theorem \ref{equilibria theorem}. Moreover, $\delta_{1}=0.1988, \delta_{2}=0.3$. By $(b)$ of Fig. 2, we know that Turing bifurcation curve $l_{1}$ intersects with Hopf bifurcation curve $H_{0}$ at $(d^{m}_2,\theta_{m}) =(0.2136,0.6627)$. Thus, numerical simulations of Hopf and steady state bifurcations are  conducted on the interval $(0,d^{m}_2)$ and $(d^{m}_2,+\infty)$, respectively. And we first carry out some numerical simulations of system \eqref{Numerical diffusion system} without diffusion at $E_{31}(0.09,0.123)$ , see Fig. \ref{localE31stable}. And Fig. \ref{diffusionE31stable} shows thats $E_{31}(0.09,0.123)$ of system \eqref{Numerical diffusion system} is asymptotically stable under some parameter conditions.
	\subsection{Hopf bifurcation and periodic solutions}
\begin{figure}[H]
		\centering
		\begin{minipage}[c]{0.5\textwidth}
			\centering
			\includegraphics[height=5cm,width=7cm]{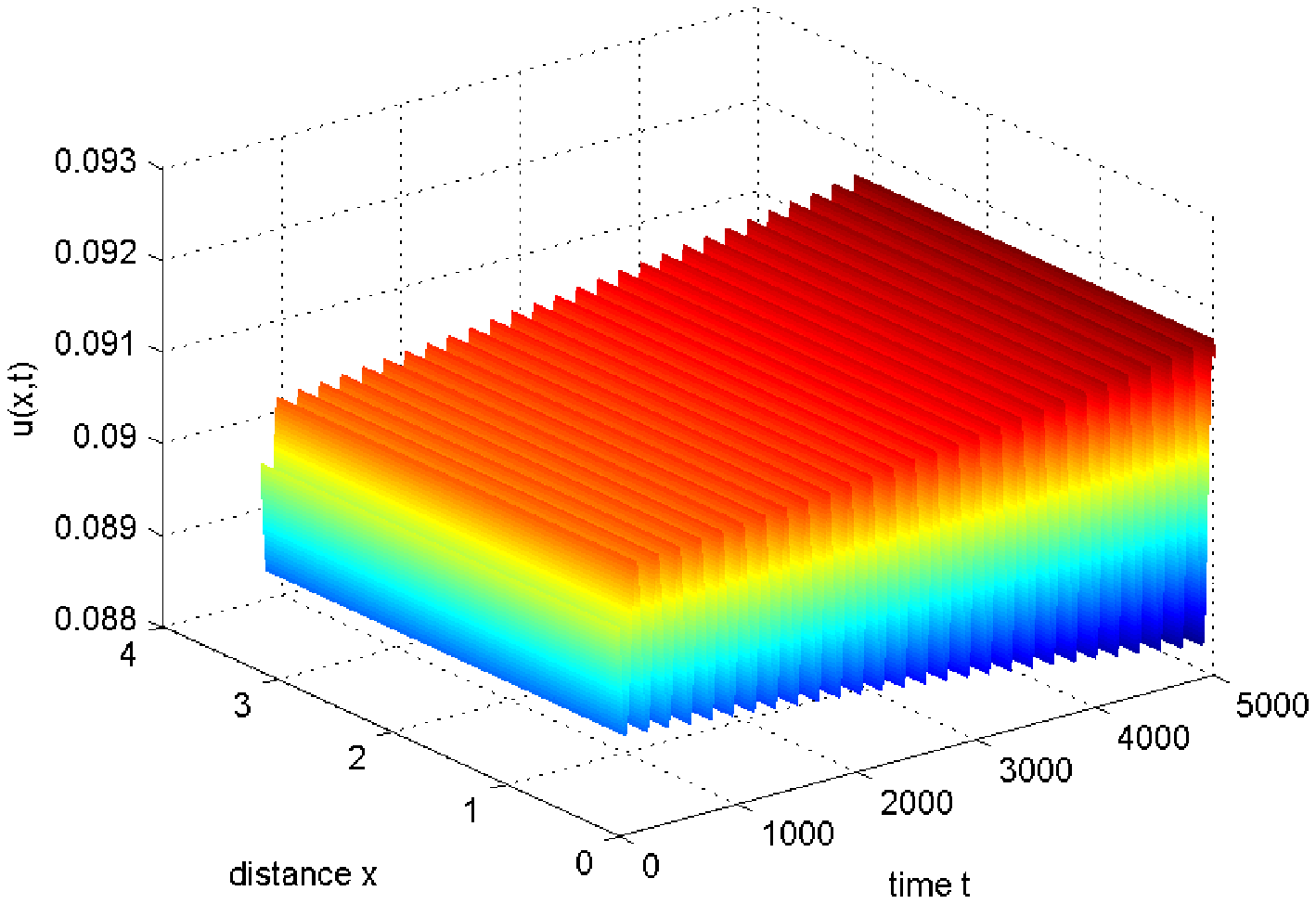}
		\end{minipage}%
		\begin{minipage}[c]{0.5\textwidth}
			\centering
			\includegraphics[height=5cm,width=7cm]{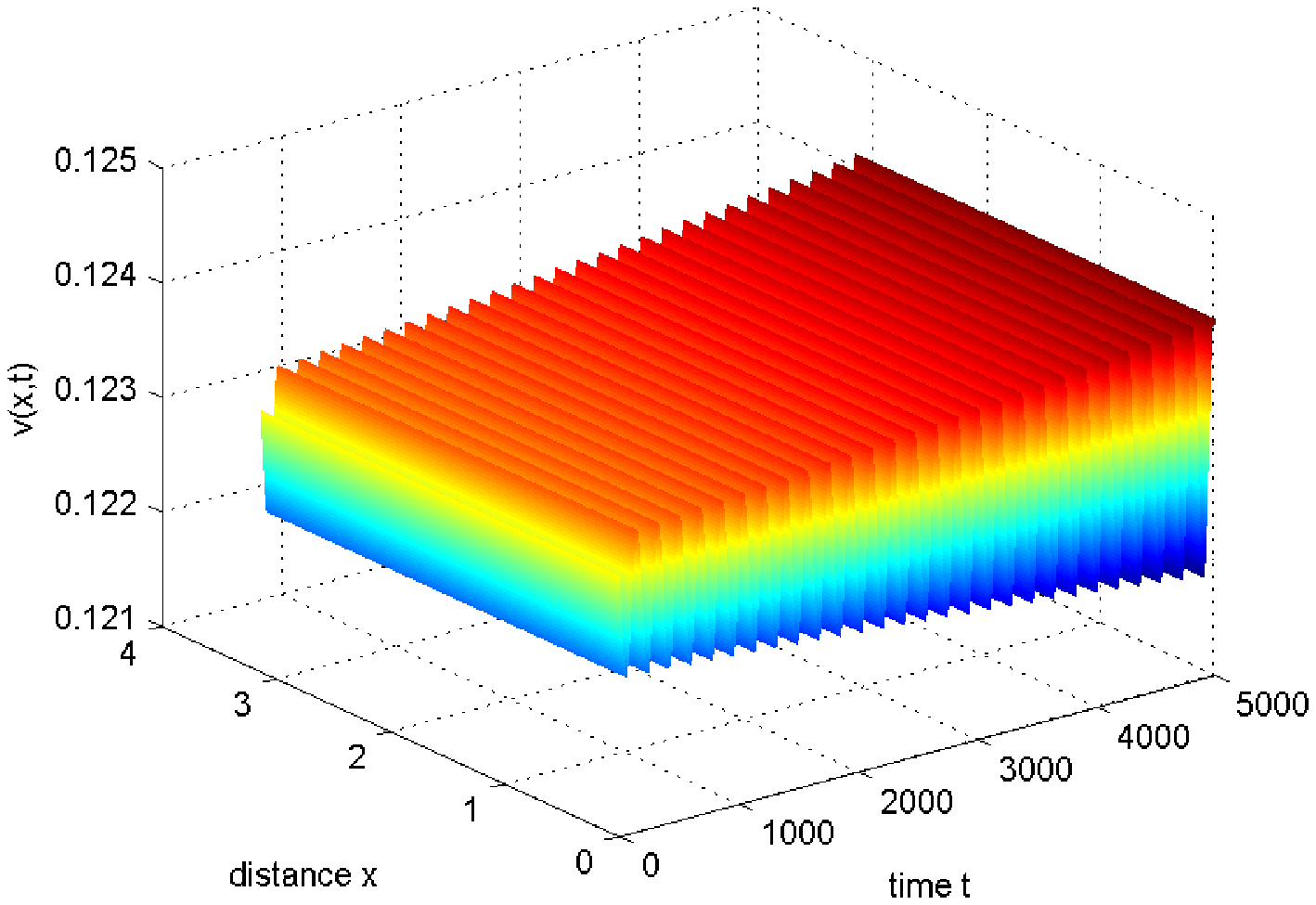}
		\end{minipage}
		\caption{For fixed parameter values
			$d_{1}=0.1,d_{2}=0.15,m=-0.5,n=\frac{50}{41},c=\frac{100}{41},\theta=0.662\textless \theta_{0}=0.6627$, a spatially homogeneous and unstable periodic solution occurs. The initial values are $(u{(x,0)},v{(x,0)})=(0.0903,0.1233)$.}\label{Hopf0uv}
	\end{figure}

	System \eqref{Numerical diffusion system} exhibits Hopf bifurcation as $\theta$ varies for $(0,d^{m}_2)$, that is, 	system \eqref{Numerical diffusion system} goes through spatially homogeneous Hopf bifurcation on $H_{0}: \theta=\theta_{0}$. To consider spatially homogeneous periodic solutions, we take $k=0, d_{2}=0.15$ and $\theta_{*}=\theta_{0}=0.6627$ in this case. Choosing $\theta=0.662$, then we have $\omega_{0}=0.035,p_{01}=1.4634+0.258\text{i},p_{02}=1,q_{01}=-1.9381 \text{i},q_{02}=0.5+2.8362\text{i}$ and
	\begin{equation*}
		f_{200} =
		\left(
		\begin{array}{c}a
			0.8734 \\ -0.7548
		\end{array}
		\right),
		f_{110}=
		\left(
		\begin{array}{c}
			-1.6667 \\ 1.1045
		\end{array}
		\right),
		f_{101} =
		\left(
		\begin{array}{c}
			0 \\  \frac{50}{41}
		\end{array}
		\right),
		f_{011} =
		\left(
		\begin{array}{c}
			0 \\   -0.3
		\end{array}
		\right),
	\end{equation*}
	\begin{equation*}
		f_{020} =
		\left(
		\begin{array}{c}
			0 \\ -3.2328
		\end{array}
		\right),
		f_{300}=
		\left(
		\begin{array}{c}
			-22.9552 \\ 12.5793
		\end{array}
		\right),
		f_{210} =
		\left(
		\begin{array}{c}
			9.2593 \\  -6.1362
		\end{array}
		\right),
		f_{120} =f_{030}=
		\left(
		\begin{array}{c}
			0 \\  0
		\end{array}
		\right).
	\end{equation*}
	Through the above analysis and calculation, the normal form takes form as
	\begin{equation}
		\dot \rho=-0.15\epsilon \rho+9.7469\rho^{3},
	\end{equation}
	which intimates $v_{01}=-0.15\textless 0$ and $v_{02}=9.7469 \textgreater 0$. By theorem \ref{Hopf bifurcation theorem}, we confirm that Hopf bifurcation on $H_{0}: \theta=\theta_{0}$ is subcritical. This indicates that an unstable spatially homogenous periodic solution arises, see Fig. \ref{Hopf0uv}.\\
	Additionally,  in order to consider spatially inhomogeneous periodic solutions, we take $k=1, d_{2}=0.002$ and $\theta_{*}=0.3227$ in this case. We choose $\theta=0.32$, then  $\omega_{1}=0.1375, p_{11}=1.4936+2.0789 \text{i}, p_{12}=1,q_{11}=-0.2405\text{i},q_{12}=0.5+0.3592\text{i}$ and
	\begin{equation*}
		f_{200} =
		\left(
		\begin{array}{c}
			0.8734 \\ -0.3675
		\end{array}
		\right),
		f_{110}=
		\left(
		\begin{array}{c}
			-1.6667 \\ 0.5379
		\end{array}
		\right),
		f_{101} =
		\left(
		\begin{array}{c}
			0 \\  \frac{50}{41}
		\end{array}
		\right),
		f_{011} =
		\left(
		\begin{array}{c}
			0 \\   -0.3
		\end{array}
		\right),
	\end{equation*}
	\begin{equation*}
		f_{020} =
		\left(
		\begin{array}{c}
			0 \\ -1.5742
		\end{array}
		\right),
		f_{300}=
		\left(
		\begin{array}{c}
			-22.9552 \\ 6.1256
		\end{array}
		\right),
		f_{210} =
		\left(
		\begin{array}{c}
			9.2593 \\  -2.9881
		\end{array}
		\right),
		f_{120} =f_{030}=
		\left(
		\begin{array}{c}
			0 \\  0
		\end{array}
		\right).
	\end{equation*}
	Therefore, the  normal form  takes form as
	\begin{equation}
		\dot \rho=-0.15\epsilon \rho+-82.6307\rho^{3},
	\end{equation}
	which intimates $v_{01}=-0.15\textless 0$ and $v_{02}=-82.6307 \textless 0$. By theorem \ref{Hopf bifurcation theorem}, we confirm that Hopf bifurcation on $H_{1}: \theta_{1}=0.3227$ is supercritical. This indicates that  a stable spatially inhomogenous periodic solution arises, see Fig. \ref{Hopf1uv}.\\
	
	\subsection{Pitchfork bifurcation and spatially inhomogeneous steady state}
\begin{figure}[H]
		\centering
		\begin{minipage}[c]{0.5\textwidth}
			\centering
			\includegraphics[height=5cm,width=7cm]{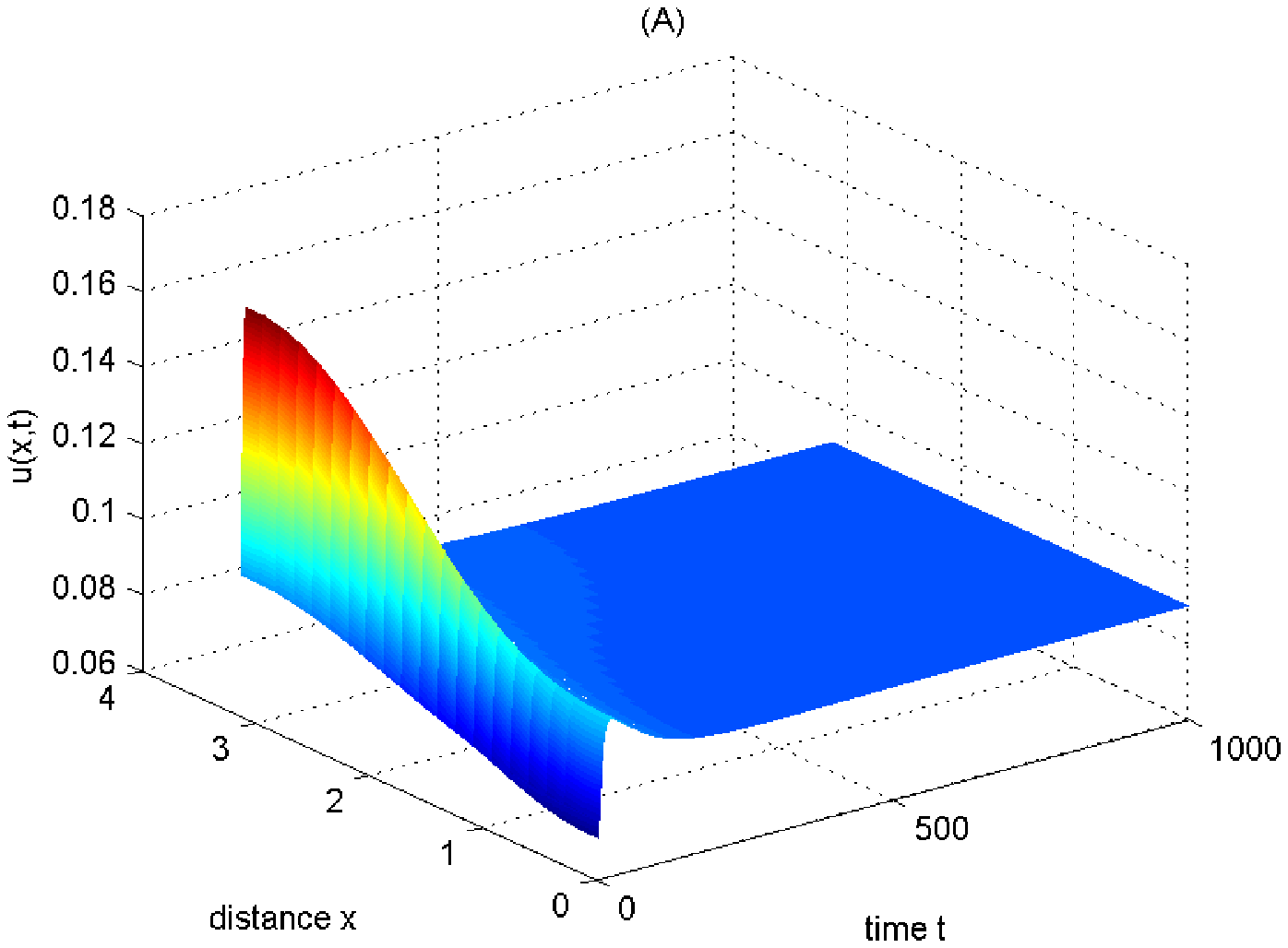}
		\end{minipage}%
		\begin{minipage}[c]{0.5\textwidth}
			\centering
			\includegraphics[height=5cm,width=7cm]{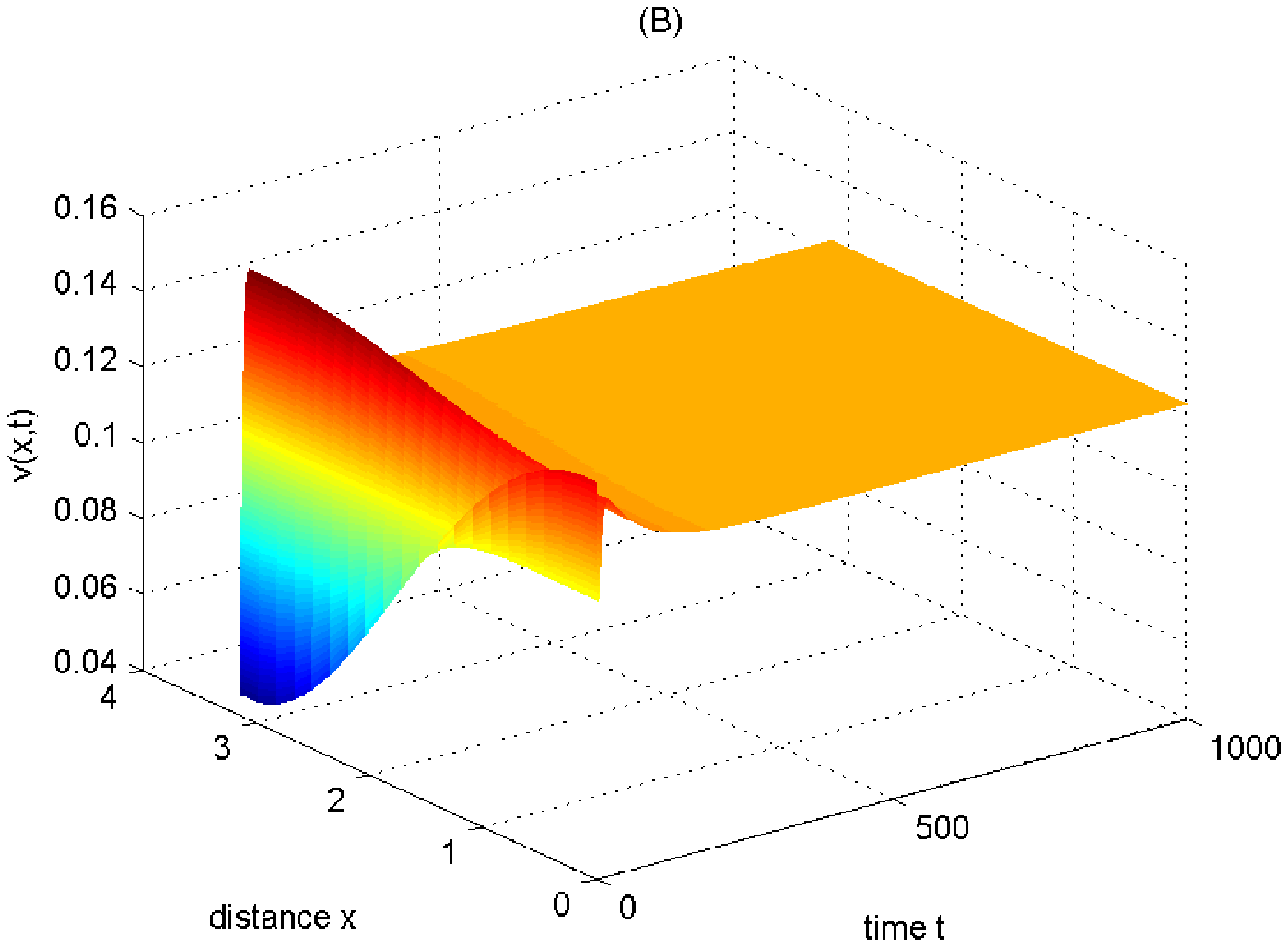}
		\end{minipage}
		\begin{minipage}[c]{0.5\textwidth}
			\centering
			\includegraphics[height=5cm,width=7cm]{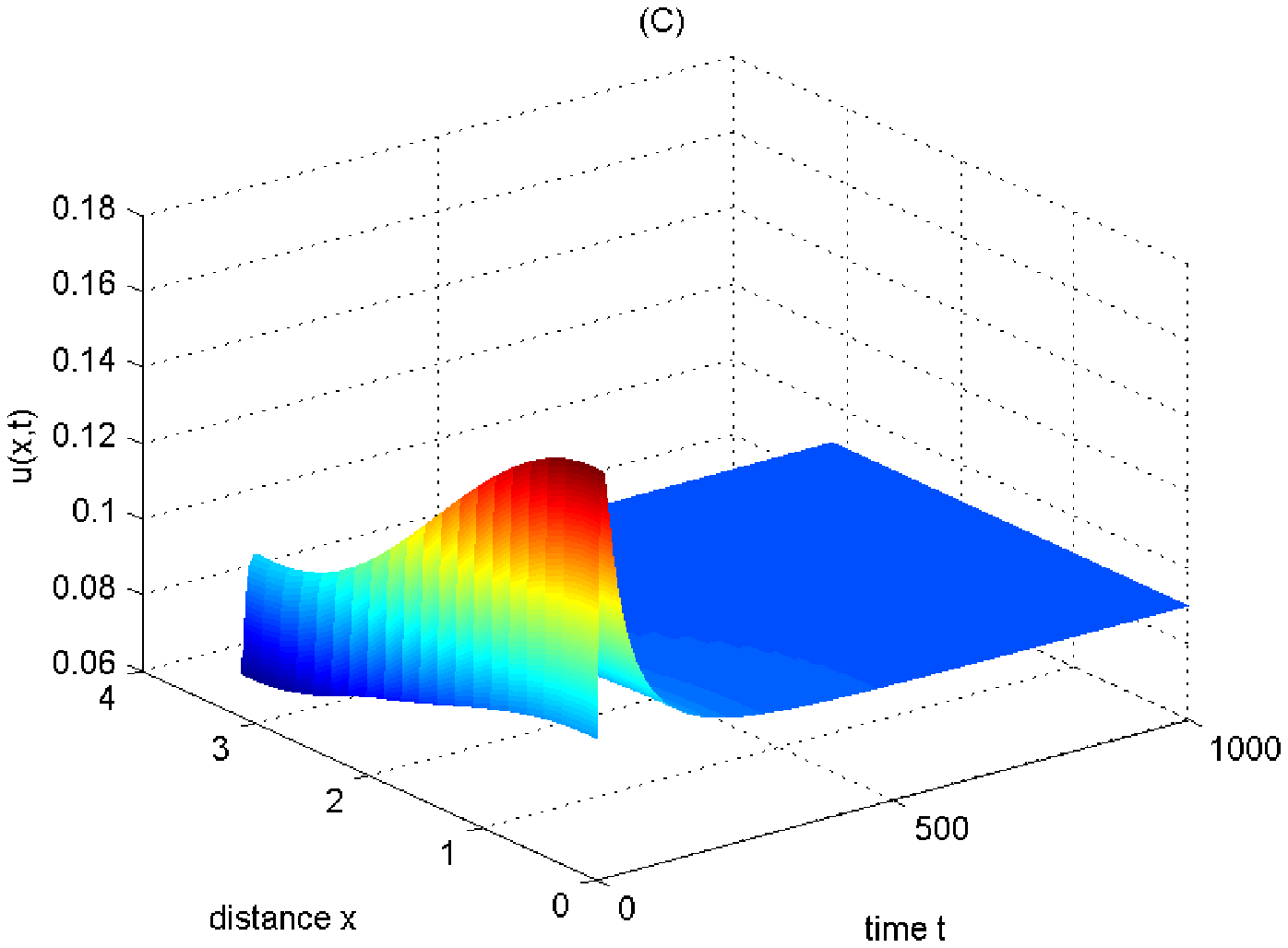}
		\end{minipage}%
		\begin{minipage}[c]{0.5\textwidth}
			\centering
			\includegraphics[height=5cm,width=7cm]{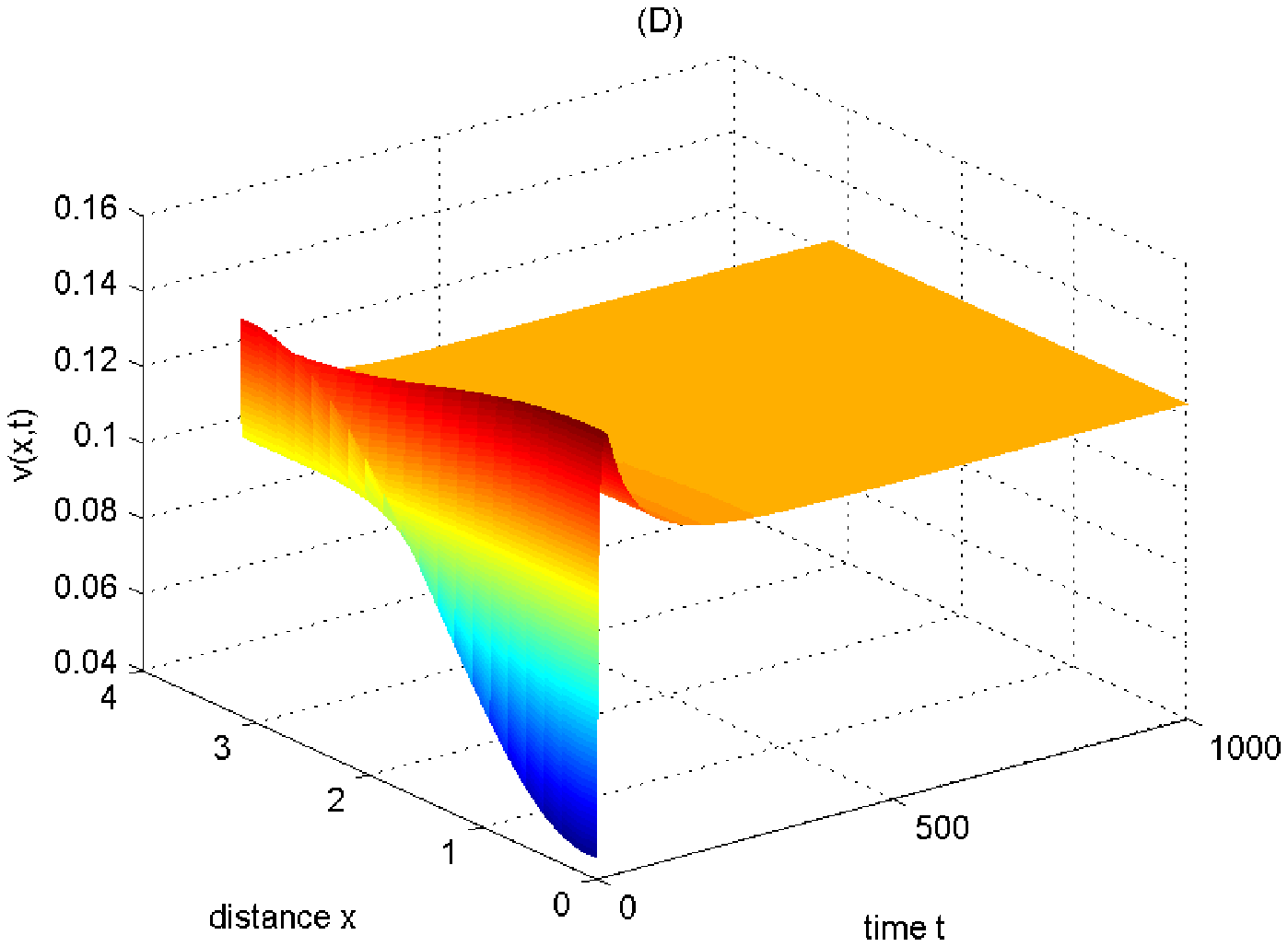}
		\end{minipage}
		\caption{ For fixed parameter values  $d_{1}=0.1, d_{2}=0.4,m=-0.5,n=\frac{50}{41},c=\frac{100}{41},\theta=1.24  \textless \theta_{T1}=1.2408$, the positive equilibrium $E_{31}(0.09,0.123)$  is unstable and two stable spatially
			inhomogeneous steady states of $\cos x$-like shape occur.  (A)-(B): The initial values are $(u{(x,0)},v{(x,0)})=(0.08+0.01 \cos x,0.1+0.1  \cos x)$. (C)-(D): The initial values are $(u{(x,0)},v{(x,0)})=(0.08-0.01 \cos x,0.1-0.1\cos x)$.}\label{Turinguv}
	\end{figure}

\begin{figure}[H]
		\centering
		\begin{minipage}[c]{0.5\textwidth}
			\centering
			\includegraphics[height=5cm,width=7cm]{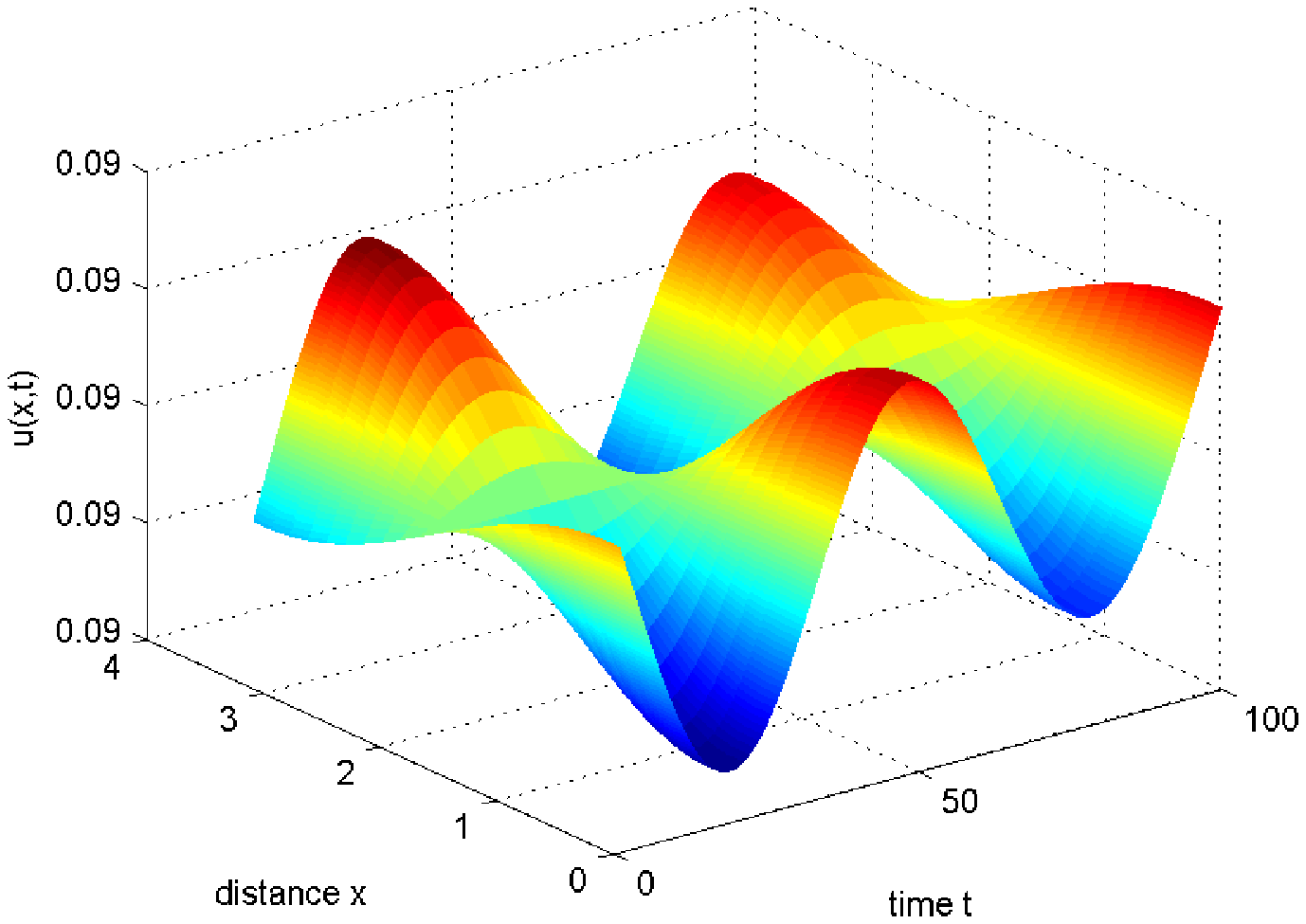}
		\end{minipage}%
		\begin{minipage}[c]{0.5\textwidth}
			\centering
			\includegraphics[height=5cm,width=7cm]{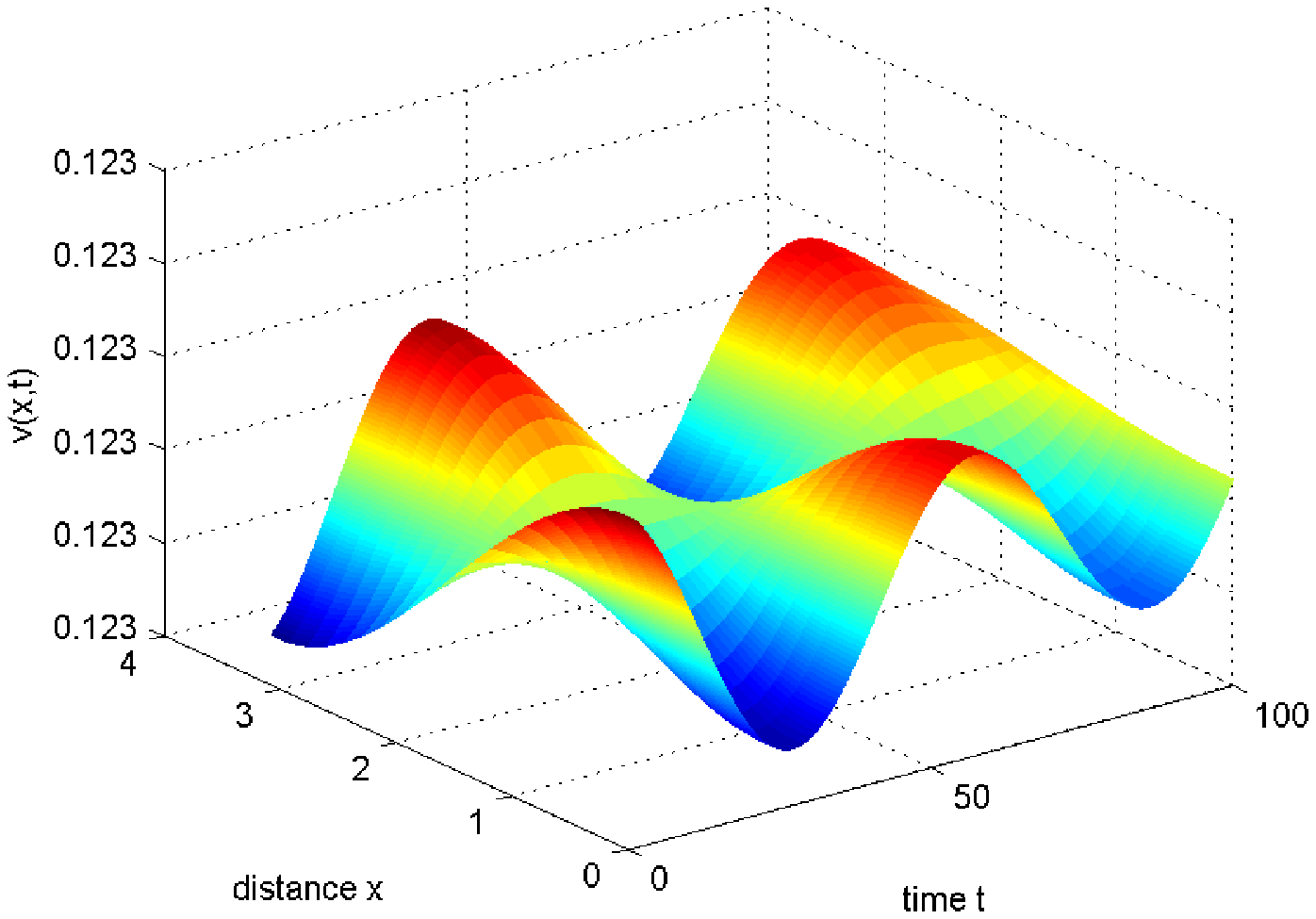}
		\end{minipage}
		\caption{For fixed parameter values
			$d_{1}=0.1,d_{2}=0.15,m=-0.5,n=\frac{50}{41},c=\frac{100}{41},\theta=0.32\textless \theta_{1}=0.3227$, a spatially inhomogeneous and stable periodic solution occurs. The initial values are $(u{(x,0)},v{(x,0)})=(0.09+0.000008 \cos x,0.123+0.000008\cos x)$.}\label{Hopf1uv}
	\end{figure}

\begin{figure}[H]
		\centering
		\begin{minipage}[c]{0.5\textwidth}
			\centering
			\includegraphics[height=5cm,width=7cm]{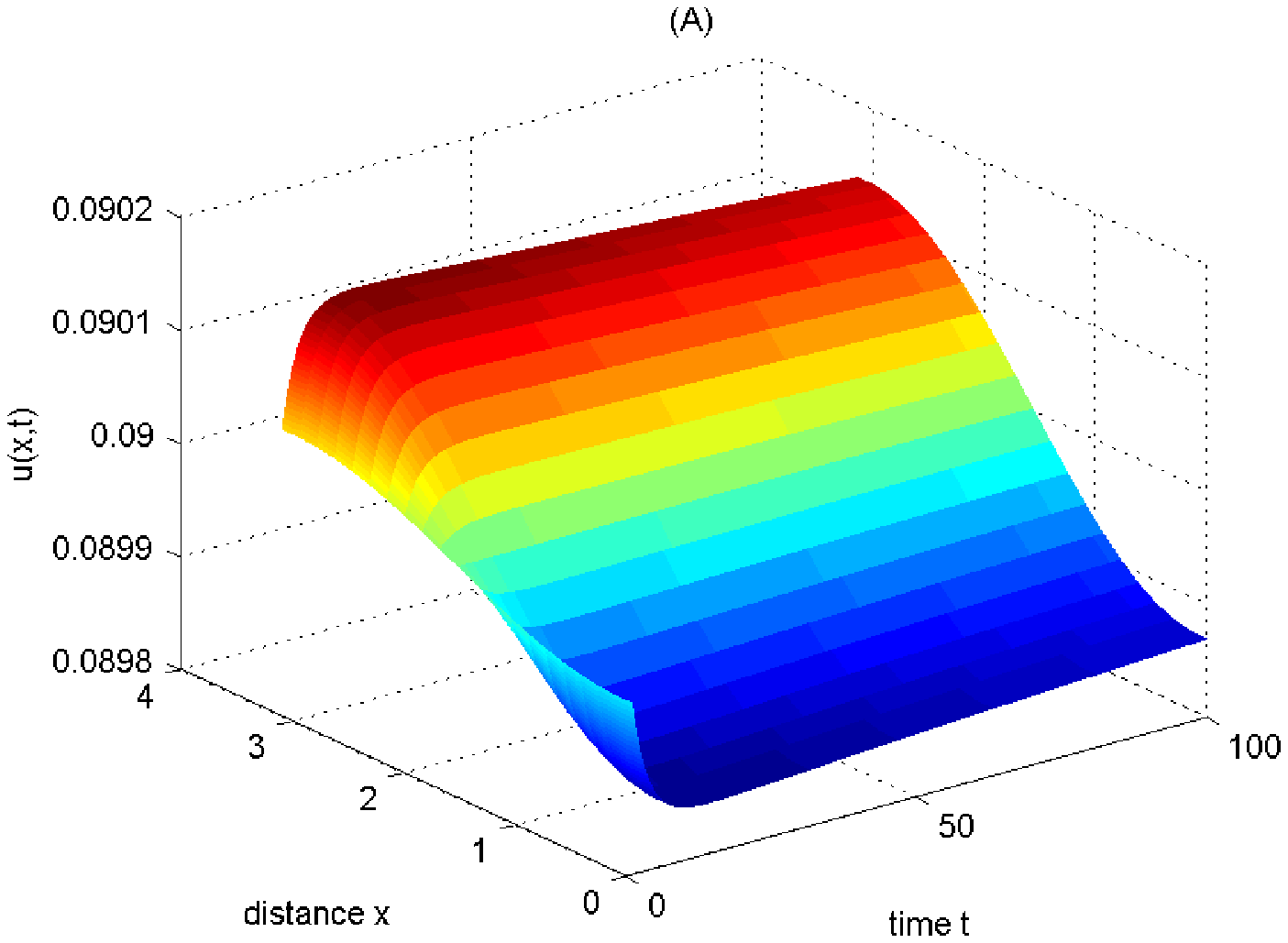}
		\end{minipage}%
		\begin{minipage}[c]{0.5\textwidth}
			\centering
			\includegraphics[height=5cm,width=7cm]{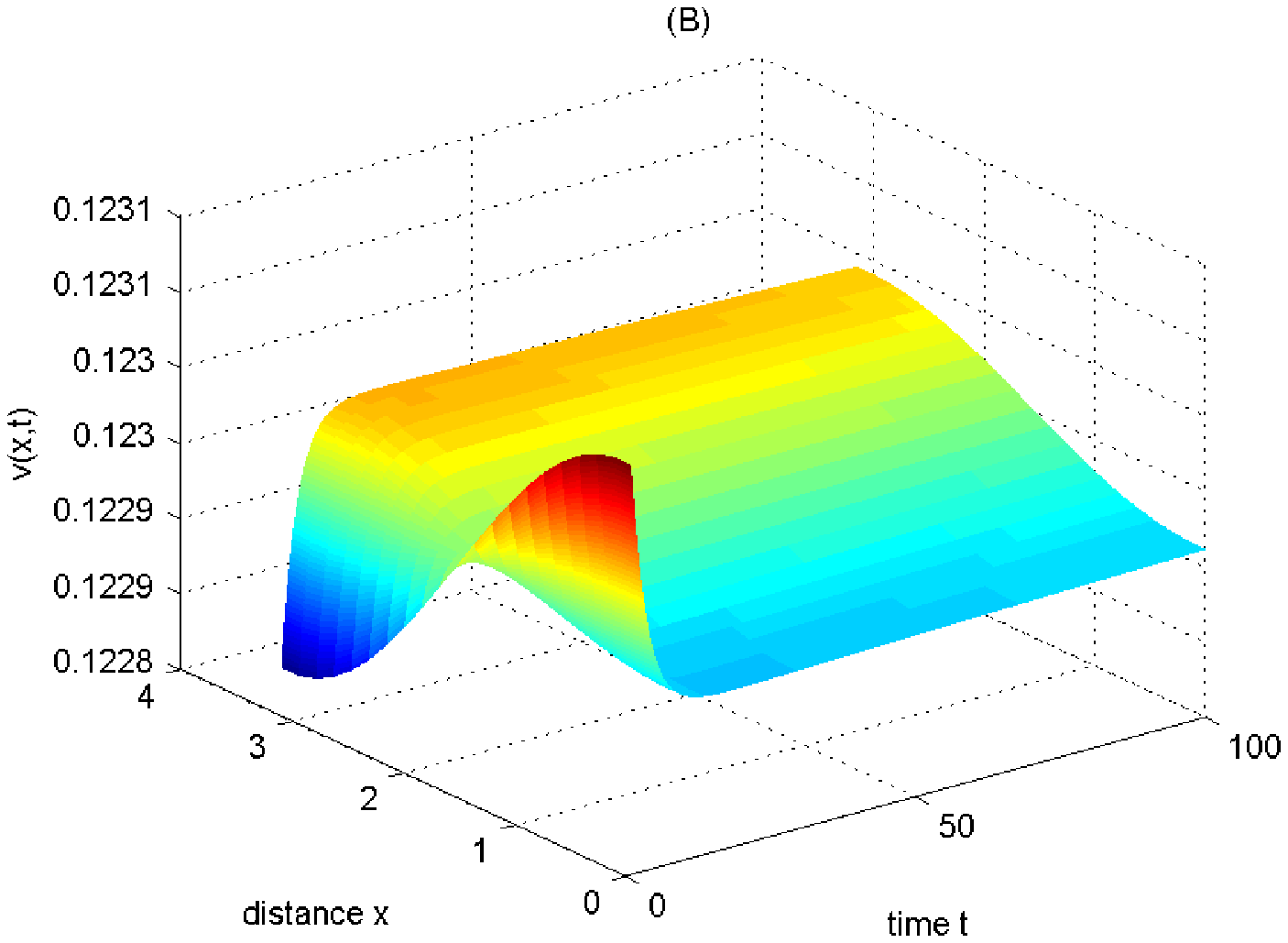}
		\end{minipage}
	\begin{minipage}[c]{0.5\textwidth}
		\centering
		\includegraphics[height=5cm,width=7cm]{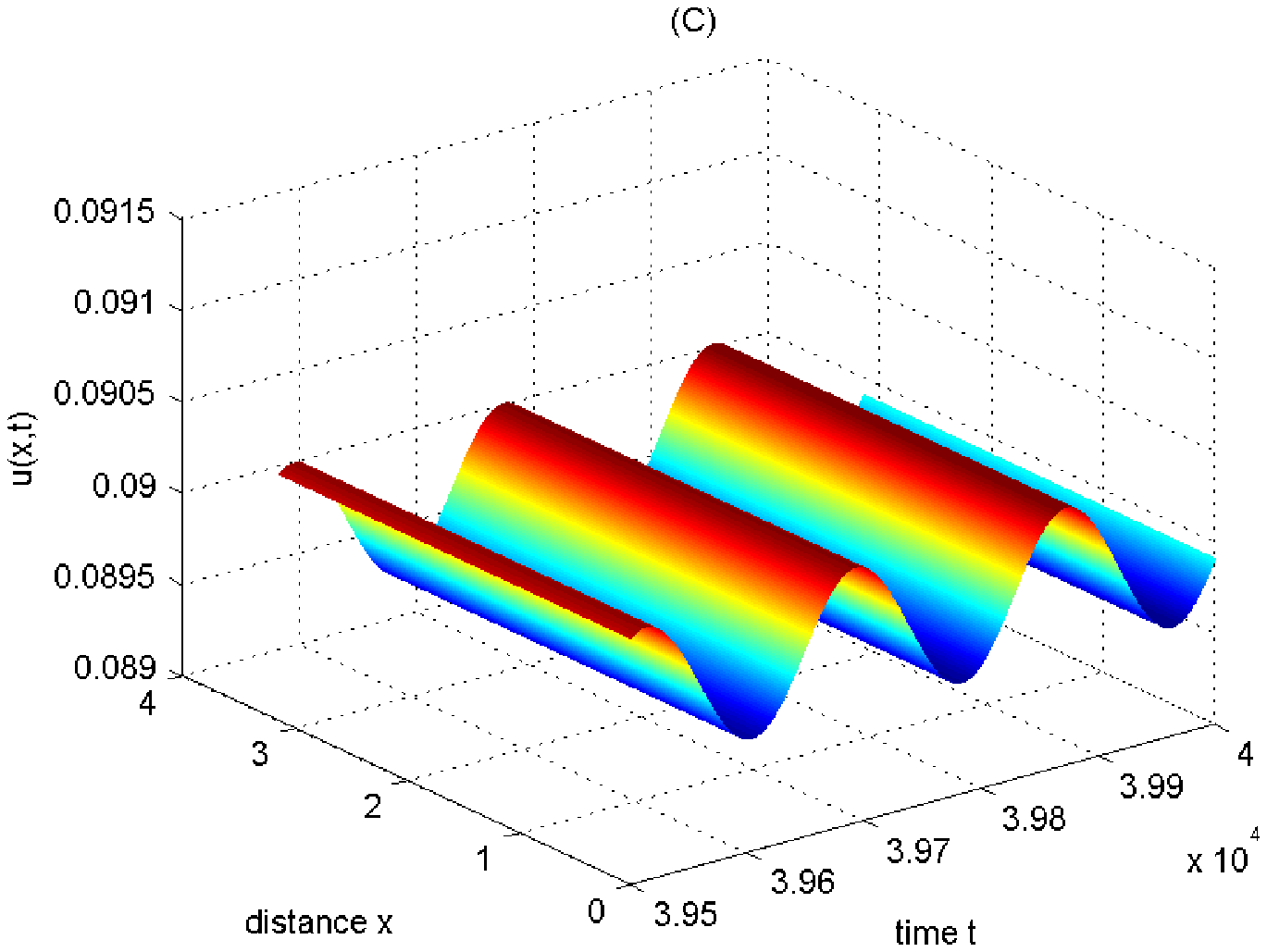}
	\end{minipage}%
	\begin{minipage}[c]{0.5\textwidth}
		\centering
		\includegraphics[height=5cm,width=7cm]{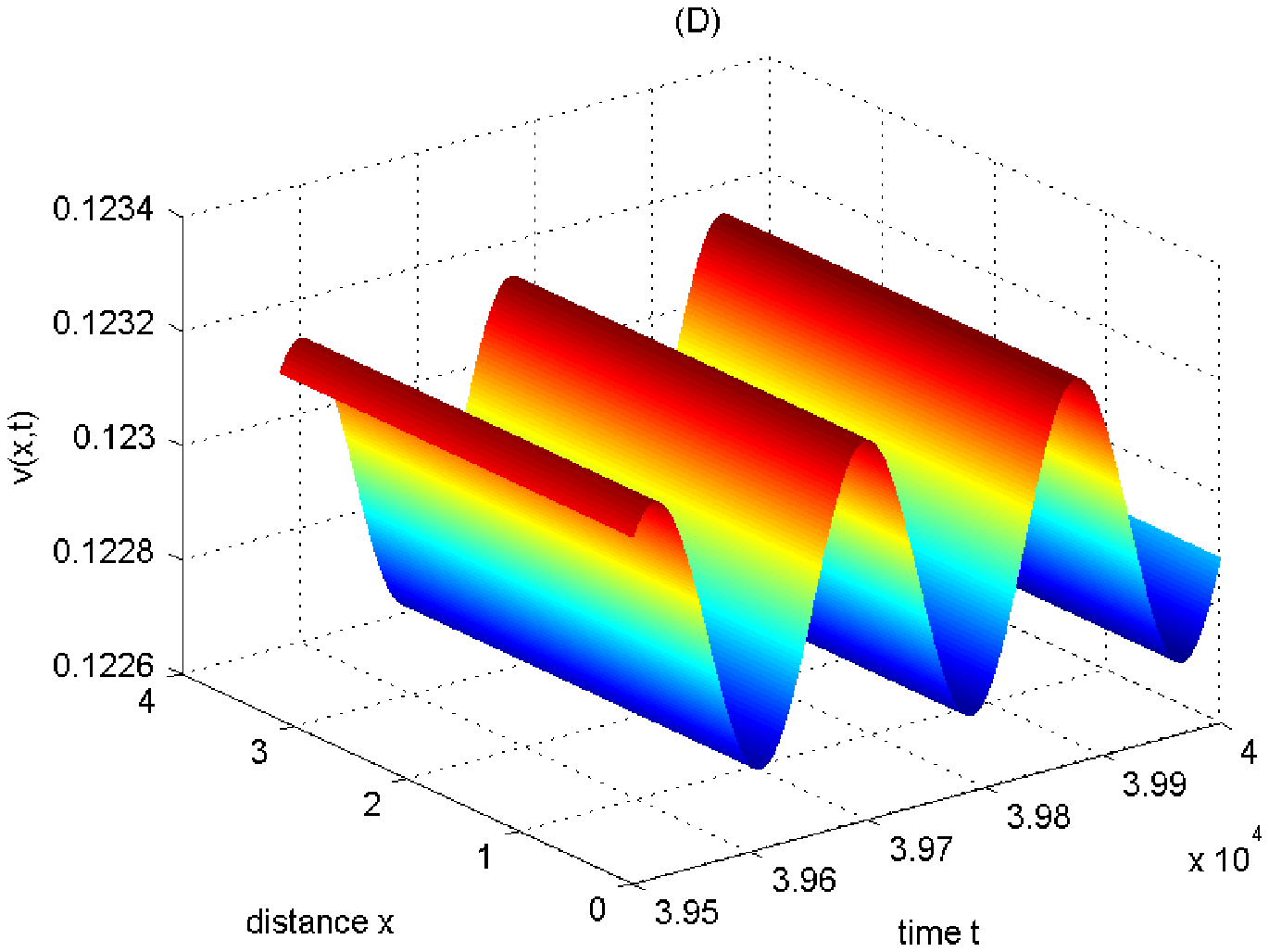}
	\end{minipage}
		\caption{For fixed parameter values  $d_{1}=0.1, d_{2}=0.23,m=-0.5,n=\frac{50}{41},c=\frac{100}{41},\theta=0.6617$, the positive equilibrium $E_{31}(0.09,0.123)$  is unstable and stable spatially
		inhomogeneous periodic solutions occur. The initial values are $(u{(x,0)},v{(x,0)})=(0.09,0.123+0.0002\cos x)$. (A)-(B) are transient
		behaviors for $u$ and $v$, respectively. (C)-(D) are long-term behaviors for $u$ and $v$, respectively.}\label{1TH}
\end{figure}

		System \eqref{Numerical diffusion system} exhibits steady state bifurcation as $\theta$ varies for $(d^{m}_2,+\infty)$. Here we take $k=1, d_{2}=0.4$ and $\theta_{*}=1.2408$ in this case. Then system \eqref{Numerical diffusion system} goes through steady state bifurcation when $\theta$ crosses the critical line $l_{1}: \theta_{T1}=1.2408$. We choose $\theta=1.24$, then  $T_{1}=-0.6732,\widetilde{p}_{11}=1,\widetilde{p}_{12}=0.3294,\widetilde{q}_{11}=1.1471,\widetilde{q}_{12}=-0.4456$ and
	\begin{equation*}
		f_{200} =
		\left(
		\begin{array}{c}
			0.8734 \\ -1.4131
		\end{array}
		\right),
		f_{110}=
		\left(
		\begin{array}{c}
			-1.6667 \\ 2.0680
		\end{array}
		\right),
		f_{101} =
		\left(
		\begin{array}{c}
			0 \\  \frac{50}{41}
		\end{array}
		\right),
		f_{011} =
		\left(
		\begin{array}{c}
			0 \\   -0.3
		\end{array}
		\right),
	\end{equation*}

	\begin{equation*}
		f_{020} =
		\left(
		\begin{array}{c}
			0 \\ -6.0526
		\end{array}
		\right),
		f_{300}=
		\left(
		\begin{array}{c}
			-22.9552 \\ 23.5518
		\end{array}
		\right),
		f_{210} =
		\left(
		\begin{array}{c}
			9.2593 \\  -11.4887
		\end{array}
		\right),
		f_{120} =f_{030}=
		\left(
		\begin{array}{c}
			0 \\  0
		\end{array}
		\right).
	\end{equation*}
Therefore, the  normal form truncated to the third-order term takes form as:
\begin{equation}
	\dot \rho=-0.4994\epsilon \rho-8.9267\times 10^{6}\rho^{3},
\end{equation}
which intimates $Q_{111}=-0.4994\textless 0$ and $Q_{130}=-8.9267\times 10^{6} \textless 0$. By theorem \ref{pitchfork bifurcation theorem}, we confirm that pitchfork bifurcation on $l_{1}: \theta_{T1}=1.2408$ is supercritical and stable. This indicates that a stable spatially inhomogeneous steady state emerges, see Fig. \ref{Turinguv}.\\
In the end,  some simulations are given around Turing-Hopf bifurcation point  $(d^{m}_2,\theta_{m}) =(0.2136,0.6627)$.  And  a family of stable spatially inhomogeneous periodic solutions  are shown in Fig. \ref{1TH}.
	
	\section{Conclusions and discussions}
	\noindent In this paper, we have investigated a  diffusive predator-prey model with multiple Allee effect, herd behavior and quadratic mortality. Moreover, we consider the quadratic mortality on predator species. We found that the dynamics of system \eqref{diffusion system} near the positive equilibria is quite rich. And we are more concerned with the spatial dynamics near the positive equilibrium $E_{31}$. The conversion rate $\theta$ of prey to predator and the diffusion rate $d_{2}$ of predator are chosen as two crucial parameters, that is, Turing instability is analyzed in the $\theta-d_{2}$ plane. And the simulations near Turing-Hopf bifurcation point are presented. We summarize our findings as follows:\\
	(1) Too large diffusion rate $d_{1}$ of prey prevents Turing instability from emerging. \\
	(2) The biomass conversion rate $\theta$ does not affect the stability of the boundary equilibria, but affects the stability of all positive equilibria except the saddle $ E_{32}$ and the occurrence of Turing instability. Thus, the biomass conversion rate is greatly important for the predator-prey system, and we can control the biological conversion rate to achieve the coexistence of the predator and prey.\\
	(3) Hopf and steady state bifurcations only occur in their respective ranges, that is,   Hopf bifurcation and steady state bifurcation are corresponding to  $d_{2}\in(0,d^{m}_2)$ and $d_{2}\in (d^{m}_2,+\infty)$ respectively. \\
	(4)  Allee effect does not alter the the local stability of the boundary equilibria $E_{1}$ and $E_{2}$ for $m\neq1$. But by considering the highest intensity of Allee effect $(m=1)$, we found that the transcritical bifurcations occur at $E_{1}$ and $E_{2}$. The strong Allee causes  $E_{31}$ to vanish, and  $E_{32}$ is always a saddle whether Allee effect is strong or weak.\\

	 \section*{Appendix 1: proof of theorem \ref{E1 theorem}-\ref{E33 theorem}}
	\noindent\textbf{Proof of theorem \ref{E1 theorem}}.\\
	The Jacobian matrix of system \eqref{local system} is
	\begin{equation}
		J(u,v)=
		\left(
		\begin{array}{cc}
			\frac{-2u^{3}+(m-3n+1)u^{2}+2n(m+1)u-mn}{(u+n)^{2}}-\frac{v}{2\sqrt{u}} & -\sqrt{u} \\
			\frac{\theta v}{2\sqrt{u}}                                              &  \theta(\sqrt{u}-2cv)
		\end{array}
		\right).
	\end{equation}
	Thus, at $E_{1}(1,0)$,
	\begin{equation}
		J(E_{1})=
		\left(
		\begin{array}{cc}
			\frac{m-1}{n+1}  &       -1 \\
			0        &     \theta
		\end{array}
		\right).
	\end{equation}
	Clearly, the eigenvalues of $J(E_{1})$ are $\lambda_{1}=\frac{m-1}{n+1}, \lambda_{2}=\theta \textgreater  0$. Hence, $E_{1}(1,0)$ is  a saddle for $-1 \textless m \textless 1$. This suggests that the strong or  weak Allee effect do not change the stability of $E_{1}$ except $m=1$.
	
	\noindent\textbf{Proof of theorem \ref{E2 theorem}}.\\
	The Jacobian matrix of system \eqref{local system} at $E_{2}(m,0)$ is
	\begin{equation}
		J(E_{2})=
		\left(
		\begin{array}{cc}
			\frac{m(1-m)}{n+m}  &       -\sqrt{m}\\
			0        &     \theta \sqrt{m}
		\end{array}
		\right).
	\end{equation}
	Obviously, the eigenvalues of $J(E_{2})$ are $\lambda_{1}=	\frac{m(1-m)}{n+m}, \lambda_{2}=\theta \sqrt{m}\textgreater 0$. Thus,  $E_{2}(m,0)$ is an unstable node for $0 \textless m\textless 1$.
	
	\noindent\textbf{Proof of theorem \ref{E30 theorem}}.\\
	 If $m_{1} \textless n \textless m_{2}$ and $(m+1)c \textgreater 1$ or $n\textless m_{1}$, then
	$\text{det}J(E_{30})=-\frac{\theta\sqrt{u_{30}}\left((1-c(m+1))u_{30}+2(mc+n)\right)}{c(u_{30}+n)}\textgreater 0$ and $\text{Tr}J(E_{30})=\frac{(3-2(m+1)c)u_{30}+4mc+5n}{2c(u_{30}+n)}-\theta\sqrt{u_{30}}$ is equivalent to  $\theta =\frac{(3-2(m+1)c)u_{30}-mc}{2c\sqrt{u_{30}}(u_{30}+n)} \triangleq \theta_{30}$.
	
	\noindent\textbf{Proof of theorem \ref{E31 theorem}}.\\
	The Jacobian matrix of system \eqref{local system} at $E_{3i}(u_{31},v_{3i})(i=0,1,2,3)$ is
	\begin{equation}
		\begin{aligned}
			J(E_{3i})&=
			\left(
			\begin{array}{cc}
				\frac{-2u^{3}_{3i}+(m-3n+1)u^{2}_{3i}+2n(m+1)u_{3i}-mn}{(u_{3i}+n)^{2}}-\frac{1}{2c} & -\sqrt{u_{3i}} \\
				\frac{\theta }{2c}                                              &  -\theta\sqrt{u_{3i}} \\
			\end{array}
			\right)  \\
			&=
			\left(
			\begin{array}{cc}
				\frac{(1-(m+1)c)u_{3i}+2(mc+n)}{c(u_{3i}+n)}+\frac{1}{2c} & -\sqrt{u_{3i}} \\
				\frac{\theta }{2c}                          &  -\theta\sqrt{u_{3i}} \\
			\end{array}
			\right).
		\end{aligned}
	\end{equation}
	Therefore,
	\begin{flalign*}
		\text{det}J(E_{3i})&=-\frac{\theta\sqrt{u_{3i}}\left((1-c(m+1))u_{3i}+2(mc+n)\right)}{c(u_{3i}+n)},\\ \text{Tr}J(E_{3i})&=\frac{(3-2(m+1)c)u_{3i}+4mc+5n}{2c(u_{3i}+n)}-\theta\sqrt{u_{3i}}.
	\end{flalign*}
	By calculation,
	$\text{det}J(E_{31})=\frac{\theta \sqrt{u_{31}}u_{31}((m+1)c-1)}{c(u_{31}-mc)}  \textgreater 0$ and $\text{Tr}J(E_{31})=0$ is equivalent to  $\theta =\frac{(3-2(m+1)c)u_{31}-mc}{2c\sqrt{u_{31}}(u_{31}-mc)} \triangleq \theta_{31}$. Hence, $E_{31}(u_{31},v_{31})$ is asymptotically stable if  $\theta \textgreater \theta_{31}$ and $E_{31}(u_{31},v_{31})$ is unstable if $\theta \textless \theta_{31}$.
	
	\noindent\textbf{Proof of theorem \ref{E32 theorem}}.\\
	If $m_{1}\textless n \textless m_{2}$, then $\text{det}J(E_{32})=
	\frac{\theta \sqrt{u_{32}}(cu^2_{32}-n-mc)}{c(u_{32}+n)}$. Clearly, $\text{det}J(E_{32})$ and $cu^2_{32}-n-mc$ have the same sign. Let $F(n)=cu^2_{32}-n-mc$, then $F^{\prime}(n)=\frac{(m+1)c-1}{\sqrt{((m+1)c-1)^2-4c(mc+n)}}-2$. Thus,  $F^{\prime}(n)\textgreater 0$ if $n \textgreater \frac{3((m+1)c-1)^{2}}{16c}-mc$ and $F^{\prime}(n)\textless 0$ if $n \textless \frac{3((m+1)c-1)^{2}}{16c}-mc$. That is,  $F(n)$ decreases monotonically on $(m_{1},\frac{3((m+1)c-1)^{2}}{16c}-mc)$ and $F(n)$ increases monotonically on $(\frac{3((m+1)c-1)^{2}}{16c}-mc,m_{2})$ with respect to $n$. This means $F(n) \textless \max \{F(m_{1}),F(m_{2})\}=0$, namely $\text{det}J(E_{32})\textless 0$.
	Cosequently, $E_{32}$ is a saddle.

	\noindent\textbf{Proof of theorem \ref{E33 theorem}}.\\
	If $n=m_{2}$, then the Jacobian matrix of system \eqref{local system} at $E_{33}(u_{33},v_{33})$ is
	\begin{equation}
		J(E_{33})=
		\left(
		\begin{array}{cc}
			\frac{(1-(m+1)c)u_{33}+2(mc+m_{2})}{c(u_{33}+m_{2})}+\frac{1}{2c} & -\sqrt{u_{33}} \\
			\frac{\theta }{2c}                          &  -\theta\sqrt{u_{33}} \\
		\end{array}
		\right)
		=
		\left(
		\begin{array}{cc}
			\frac{1}{2c}           &      -\sqrt{u_{33}} \\
			\frac{\theta }{2c}     &  -\theta\sqrt{u_{33}} \\
		\end{array}
		\right).
	\end{equation}
	Hence, $\text{det}J(E_{33})=0$.
	And $\text{Tr}J(E_{33})=0$ is equivalent to  $\theta =\frac{1}{2c\sqrt{u_{33}}} \triangleq \theta_{33}$.
	Obviously,  one of the eigenvalues of $J(E_{33})$ is 0 for $\theta \neq \theta_{33}$. To investigate the stability of $E_{33}$, we first shift $(u_{33},v_{33})$ to the origin by the transformation $(\widehat{u},\widehat{v})=(u-u_{33},v-v_{33})$ and perform Taylor expansion of system \eqref{local system} at the origin to the third order. After dropping the bars, then we get
	\begin{equation}\label{gongshi4}
		\begin{cases}
			\frac{du}{dt}=a_{1}u+a_{2}v+a_{3}u^{2}+a_{4}uv+a_{5}u^{2}v+a_{6}u^{3}+O({|u,v|^{4}}), \\
			\frac{dv}{dt}=b_{1}u+b_{2}v+b_{3}u^{2}+b_{4}uv+b_{5}v^{2}+b_{6}u^{3}+b_{7}u^{2}v+O({|u,v|^{4}}),
		\end{cases}
	\end{equation}
	where
	\begin{equation*}
		\begin{aligned}
			a_{1}&=\frac{1}{2c}, a_{2}= -\sqrt{u_{33}}, a_{3}=\frac{-u^{3}_{33}-3nu^{2}_{33}-3n^{2}u_{33}+(m+1)n^{2}+mn}{(u_{33}+n)^{3}}+\frac{1}{8cu_{33}},\\
			a_{4}&=-\frac{1}{2\sqrt{u_{33}}},
			a_{5}=\frac{1}{8u^{\frac{3}{2}}_{33}},
			b_{1}=\frac{\theta}{2c},
			b_{2}=-\theta\sqrt{u_{33}},
			b_{3}=-\frac{\theta}{8cu_{33}},
			b_{4}=\frac{\theta}{2\sqrt{u_{33}}},\\
			b_{5}&=-c\theta,
			b_{6}=\frac{\theta}{16cu^{2}_{33}},
			b_{7}=-\frac{\theta}{8u^{\frac{3}{2}}_{33}}.
		\end{aligned}
	\end{equation*}
	Applying the transformation
	\begin{equation*}
		\left(
		\begin{array}{c}
			\widehat{u}  \\ \widehat{v}
		\end{array}
		\right)
		=
		\left(
		\begin{array}{cc}
			1& -\frac{1}{\theta}\\
			0& 1
		\end{array}
		\right)
		\left(
		\begin{array}{c}
			u\\ v
		\end{array}
		\right)
	\end{equation*}
	and dropping the bars , then \eqref{gongshi4} is transformed to the following form:
	\begin{equation}\label{gongshi5}
		\begin{cases}
			\frac{du}{dt}=c_{1}u^{2}+c_{2}uv+c_{3}v^{2}+c_{4}u^{2}v+c_{5}uv^{2}+c_{6}u^{3}+c_{7}v^{3}+O({|u,v|^{4}}), \\
			\frac{dv}{dt}=d_{1}u+d_{2}v+d_{3}u^{2}+d_{4}v^{2}+d_{5}uv+d_{6}u^{2}v+d_{7}uv^{2}+d_{8}u^{3}+d_{9}v^{3}+O({|u,v|^{4}}),
		\end{cases}
	\end{equation}
	where
	\begin{equation*}
		\begin{aligned}
			c_{1}&=a_{3}-\frac{b_{3}}{\theta},  c_{2}=a_{4}-\frac{b_{4}}{\theta}+\frac{2}{\theta}(a_{3}-\frac{b_{3}}{\theta}),
			c_{3}=\frac{1}{\theta^{2}}(a_{3}-\frac{b_{3}}{\theta})+\frac{1}{\theta}(a_{4}-\frac{b_{4}}{\theta})-\frac{b_{5}}{\theta},\\
			c_{4}&=a_{5}-\frac{b_{7}}{\theta}+\frac{3}{\theta}(a_{6}-\frac{b_{6}}{\theta}),
			c_{5}=\frac{2}{\theta}(a_{5}-\frac{b_{7}}{\theta})+\frac{3}{\theta^{2}}(a_{6}-\frac{b_{6}}{\theta}),
			c_{6}=a_{6}-\frac{b_{6}}{\theta},\\
			c_{7}&=\frac{1}{\theta^{2}}(a_{5}-\frac{b_{7}}{\theta})+\frac{1}{\theta^{3}}(a_{6}-\frac{b_{6}}{\theta}),
			d_{1}=b_{1},
			d_{2}=b_{2}+\frac{b_{1}}{\theta},
			d_{3}=b_{3},
			d_{4}=b_{4}+b_{5}+\frac{2b_{3}}{\theta},\\
			d_{5}&=b_{4}+\frac{2b_{3}}{\theta},
			d_{6}=b_{7}+\frac{3b_{6}}{\theta},
			d_{7}=\frac{3b_{6}}{\theta^{2}}+\frac{2b_{7}}{\theta},
			d_{8}=b_{6},
			d_{9}=\frac{b_{6}}{\theta^{3}}+\frac{b_{7}}{\theta^{2}}.
		\end{aligned}
	\end{equation*}
	Since 	$\theta \neq \theta_{33}$ implies $d_{2} \neq 0$,  we apply the transformation
	\begin{equation*}
		\left(
		\begin{array}{c}
			\widehat{u}\\ \widehat{v}
		\end{array}
		\right)
		=
		\left(
		\begin{array}{cc}
			1& 0\\ d_{1}& d_{2}
		\end{array}
		\right)
		\left(
		\begin{array}{c}
			u \\ v
		\end{array}
		\right)
	\end{equation*}
	and dropping the bars, then \eqref{gongshi5} is transformed  to the following form:
	\begin{equation}\label{gongshi6}
		\begin{cases}
			\frac{du}{dt}=k_{1}u^{2}+k_{2}uv+k_{3}v^{2}+k_{4}u^{2}v+k_{5}uv^{2}+k_{6}u^{3}+k_{7}v^{3}+O({|u,v|^{4}}), \\
			\frac{dv}{dt}=v+h_{1}u^{2}+h_{2}uv+h_{3}v^{2}+h_{4}u^{2}v+h_{5}uv^{2}+h_{6}u^{3}+h_{7}v^{3}+O({|u,v|^{4}}),
		\end{cases}
	\end{equation}
	where
	\begin{equation*}
		\begin{aligned}
			k_{1}&=c_{1}-\frac{c_{2}d_{1}}{d_{2}}+\frac{c_{3}d^{2}_{1}}{d^{2}_{2}},
			k_{2}=\frac{c_{2}}{d_{2}}-\frac{2c_{3}d_{1}}{d^{2}_{2}},
			k_{3}=\frac{c_{3}}{d^{2}_{2}},
			k_{4}=\frac{c_{4}}{d_{2}}-\frac{2c_{5}d_{1}}{d^{2}_{2}}+\frac{3c_{7}d^{2}_{1}}{d^{3}_{2}},\\
			k_{5}&=\frac{c_{5}}{d^{2}_{2}}-\frac{3c_{7}d_{1}}{d^{3}_{2}},
			k_{6}=c_{6}-\frac{c_{4}d_{1}}{d_{2}}+\frac{c_{5}d^{2}_{1}}{d^{2}_{2}}-\frac{c_{7}d^{3}_{1}}{d^{3}_{2}},
			k_{7}=\frac{c_{7}}{d^{3}_{2}},\\
			h_{1}&=d_{3}+\frac{d_{4}d^{2}_{1}}{d^{2}_{2}}-\frac{d_{1}d_{5}}{d_{2}},
			h_{2}=\frac{d_{5}}{d_{2}}-\frac{2d_{1}d_{4}}{d^{2}_{2}},
			h_{3}=\frac{d_{4}}{d^{2}_{2}},
			h_{4}=\frac{d_{6}}{d_{2}}+\frac{3d^{2}_{1}d_{9}}{d^{3}_{2}}-\frac{2d_{1}d_{7}}{d^{2}_{2}},\\
			h_{5}&=\frac{d_{7}}{d^{2}_{2}}-\frac{3d_{1}d_{9}}{d^{3}_2},
			h_{6}=\frac{d_{7}d^{2}_{1}}{d^{2}_{2}}-\frac{d_{1}d_{6}}{d_{2}}+d_{8}-\frac{d^{3}_{1}d_{9}}{d^{3}_2},
			h_{7}=\frac{d_{9}}{d^{3}_2}.
		\end{aligned}
	\end{equation*}
	By $\frac{dv}{dt}=0$, we get the implicit function $v=-h_{1}u^{2}+(h_{1}h_{2}-h_{6})u^{3}+\cdots$. Then  $\frac{du}{dt}=k_{1}u^{2}-k_{2}h_{1}u^{3}+\cdots$. If $k_{1}\neq 0$, by Theorem 7.1  in Zhang et al. \cite{zhangzhifen} [Page 131], then $E_{33}$ is a saddle node.\\
	\indent Next we prove the  case $(ii)$. If  $\theta = \theta_{33}$,  then  the eigenvalues of $J(E_{33})$ are both 0.  Applying the transformation 	
	\begin{equation*}
		\left(
		\begin{array}{c}
			\widehat{u}  \\ \widehat{v}
		\end{array}
		\right)
		=
		\left(
		\begin{array}{cc}
			1      &       0\\
			-\theta &       1
		\end{array}
		\right)
		\left(
		\begin{array}{c}
			u\\ v
		\end{array}
		\right)
	\end{equation*}
	to transform \eqref{gongshi4} to the following form after dropping the  bars.
	\begin{equation}\label{gongshi7}
		\begin{cases}
			\frac{du}{dt}=\alpha_{1}v+\alpha_{2}u^{2}+\alpha_{3}uv+\alpha_{4}u^{2}v+\alpha_{5}u^{3}+O({|u,v|^{4}}), \\
			\frac{dv}{dt}=\beta_{1}u^{2}+\beta_{2}uv+\beta_{3}v^{2}+\beta_{4}u^{2}v+\beta_{5}u^{3}+O({|u,v|^{4}}),
		\end{cases}
	\end{equation}
	where
	\begin{equation*}
		\begin{aligned}
			\alpha_{1}&=a_{2},\alpha_{2}=a_{3}+a_{4}\theta,\alpha_{3}=a_{4},\alpha_{4}=a_{5},\alpha_{5}=a_{6}+a_{5}\theta,\\
			\beta_{1}&=b_{3}+b_{4}\theta+b_{5}\theta^{2}-a_{3}\theta-a_{4}\theta^{2},\beta_{2}=b_{4}+2b_{5}\theta-a_{4}\theta,\beta_{3}=b_{5},\beta_{4}=b_{7}-a_{5}\theta,\\
			\beta_{5}&=b_{6}+b_{7}\theta-a_{5}\theta^{2}-a_{6}\theta.
		\end{aligned}
	\end{equation*}
	Applying  time rescaling $\tau =\alpha_{1}t$, system \eqref{gongshi7} is transformed into the following form:
	\begin{equation}\label{gongshi8}
		\begin{cases}
			\begin{aligned}
				\frac{du}{d\tau}&=v+\frac{\alpha_{2}}{\alpha_{1}}u^{2}+\frac{\alpha_{3}}{\alpha_{1}}uv+\frac{\alpha_{4}}{\alpha_{1}}u^{2}v+\frac{\alpha_{5}}{\alpha_{1}}u^{3}+O({|u,v|^{4}}) \\
				&\triangleq  P(u,v)       ,  \\
				\frac{dv}{d\tau}&=\frac{\beta_{1}}{\alpha_{1}}u^{2}+\frac{\beta_{2}}{\alpha_{1}}uv+\frac{\beta_{3}}{\alpha_{1}}v^{2}+\frac{\beta_{4}}{\alpha_{1}}u^{2}v+\frac{\beta_{5}}{\alpha_{1}}u^{3}+O({|u,v|^{4}})\\
				&\triangleq Q(u,v),
			\end{aligned}
		\end{cases}
	\end{equation}
	From $P(u,v)=0$, we get the implicit function $I(u)=-\frac{\alpha_{2}}{\alpha_{1}}u^{2}+(\frac{\alpha_{2}\alpha_{3}}{\alpha^{2}_{1}}-\frac{\alpha_{5}}{\alpha_{1}})u^{3}+\cdots$, then $Q(u,I(u))=\frac{\beta_{1}}{\alpha_{1}}u^{2}+(\frac{\beta_{5}}{\alpha_{1}}-\frac{\alpha_{2}\beta_{2}}{\alpha^{2}_{1}})u^{3}+\cdots$. Thus, $\frac{\partial P(u,I(u))}{\partial u}+\frac{\partial Q(u,I(u))}{\partial v}=\frac{2\alpha_{2}}{\alpha_{1}}u+\cdots$. From theorem 7.3  and its corollary in Zhang et al. \cite{zhangzhifen} [Pages 152-155],  we have $k=2M=2, M=1, a_{k}=\frac{\beta_{1}}{\alpha_{1}},N=1,B_{N}=\frac{2\alpha_{2}}{\alpha_{1}}$. Thus,   $E_{33}$ is a degenerated singularity if $\alpha_{2}=0$, $E_{33}$ is a  saddle node if $\alpha_{2} \neq0$.\\
	Now taking
	\begin{equation}\label{gongshi9}
		\begin{cases}
			\frac{d \widetilde{u}}{d\tau}= u,       \\
			\frac{d \widetilde{v}}{d\tau}=\frac{\alpha_{2}}{\alpha_{1}}u^{2}+\frac{\alpha_{3}}{\alpha_{1}}uv+\frac{\alpha_{4}}{\alpha_{1}}u^{2}v+\frac{\alpha_{5}}{\alpha_{1}}u^{3}+O({|u,v|^{4}}),
		\end{cases}
	\end{equation}
	and dropping the bars, then we get
	\begin{equation}\label{gongshi10}
		\begin{cases}
			\frac{du}{d\tau}= v,       \\
			\frac{dv}{d\tau}=\frac{\beta_{1}}{\alpha_{1}}u^{2}+(\frac{\alpha_{3}}{\alpha_{1}}+\frac{2\alpha_{2}}{\alpha_{1}})uv+O({|u,v|^{2}}).
		\end{cases}
	\end{equation}
	Using theorem 3 in Perko \cite{Perko}, $E_{33}$ is a cusp of codimension at least 3 if $(\frac{\alpha_{3}}{\alpha_{1}}+\frac{2\alpha_{2}}{\alpha_{1}})$ i.e. $2\alpha_{2}+\alpha_{3}=0$ and $E_{33}$ is a cusp of codimension 2 if $(\frac{\alpha_{3}}{\alpha_{1}}+\frac{2\alpha_{2}}{\alpha_{1}})\neq0$ i.e. $2\alpha_{2}+\alpha_{3} \neq0$.		
	
	\section*{Appendix 2: proof of theorem \ref{Transcritical theorem1}-\ref{Saddle-node theorem}}
	\noindent\textbf{Proof of theorem \ref{Transcritical theorem1}-\ref{Transcritical theorem2}}.\\
  Apparently, one of the eigenvalues of $J(E_{1})$ is 0 for $m =1$. Set $V$ and $W$ are  the eigenvectors of $J(E_{1})$ and $J(E_{1})^{T}$  corresponding  the zero eigenvalue, respectively. Then $V$ and $W$ are as follows:
	\begin{equation*}
		V =
		\left(
		\begin{array}{c}
			V_{1} \\V_{2}
		\end{array}
		\right)
		=
		\left(
		\begin{array}{c}
			1 \\ 0
		\end{array}
		\right);
		W=
		\left(
		\begin{array}{c}
			W_{1} \\W_{2}
		\end{array}
		\right)
		=
		\left(
		\begin{array}{c}
			\theta \\ 1
		\end{array}
		\right).
	\end{equation*}
	Furthermore, we get
	\begin{equation*}
		G_{m}(E_{1},m_{TC})=
		\left(
		\begin{array}{c}
			\frac{u(1-u)}{u+n} \\ 0
		\end{array}
		\right)
		_{(E_{1},m_{TC})}
		=\left(
		\begin{array}{c}
			0 \\ 0
		\end{array}
		\right),
	\end{equation*}
	\begin{equation*}
		DG_{m}(E_{1},m_{TC})V=
		\left(
		\begin{array}{cc}
			-\frac{u^{2}+2nu-n}{(u+n)^{2}} &  0 \\
			0 & 0 \\
		\end{array}
		\right)
		\left(
		\begin{array}{c}
			1 \\ 0
		\end{array}
		\right)
		_{(E_{1},m_{TC})}
		=\left(
		\begin{array}{c}
			-\frac{1}{n+1}  \\ 0
		\end{array}
		\right),
	\end{equation*}
	\begin{equation*}
		\begin{aligned}
			D^{2}G(E_{1},m_{TC})(V,V)&=
			\left(
			\begin{array}{c}
				\frac{{\partial }^{2}{G}_{1}}{\partial {x}^{2}}V_{1}^{2}+\frac{{\partial }^{2}{G}_{1}}{\partial {x} \partial {y}}V_{1}V_{2}+\frac{{\partial }^{2}{G}_{1}}{\partial {y}^{2}}V_{2}^{2}  \\
				\frac{{\partial }^{2}{G}_{2}}{\partial {x}^{2}}V_{1}^{2}+\frac{{\partial }^{2}{G}_{2}}{\partial {x} \partial {y}}V_{1}V_{2}+\frac{{\partial }^{2}{G}_{2}}{\partial {y}^{2}}V_{2}^{2}
			\end{array}
			\right)
			_{(E_{1},m_{TC})}                          \\
			&= \left(
			\begin{array}{c}
				-\frac{2}{n+1}  \\  0
			\end{array}
			\right).                                 \\
		\end{aligned}
	\end{equation*}
	It is apparent that V and W satisfy
	\begin{equation*}
		W^{T}F_{c}(E_{1},c_{TC})=0,
	\end{equation*}
	\begin{equation*}
		W^{T}[DF_{c}(E_{1},c_{TC})]=-\frac{\theta}{n+1}\neq0,
	\end{equation*}
	\begin{equation*}
		W^{T}[D^{2}F(E_{1},c_{TC})(V,V)]=-\frac{2\theta}{n+1}\neq0.
	\end{equation*}
	From Sotomayor's theorem in Perko \cite{Perko}, the transcritical  bifurcation  appears at $E_{1}(1,0)$. \\
	If $m=1$, then $E_{2}(m,0)$ coincides with $E_{1}(1,0)$.
	Analogously, we conclude that the transcritical bifurcation occurs at $E_{2}(m,0)$.\\
   \noindent\textbf{Proof of theorem \ref{Saddle-node theorem}}.\\
	Apparently, one of the eigenvalues of $J(E_{33})$ is 0 for $\theta\neq \theta_{33}$. \\
	And let $\kappa=\frac{2u^{3}_{33}+6n_{SN}u^{2}_{33}+6n^{2}_{SN}u_{33}-(2m+2)n^{2}_{SN}-2mn_{SN}}{(u_{33}+n_{SN})^{3}}-\frac{1}{cu_{33}}$. Set $V$ and $W$ are  the eigenvectors of $J(E_{21})$ and $J(E_{21})^{T}$  corresponding  the zero eigenvalue, respectively. Then $V$ and $W$ are as follows:
	\begin{equation*}
		V =
		\left(
		\begin{array}{c}
			V_{1} \\V_{2}
		\end{array}
		\right)
		=
		\left(
		\begin{array}{c}
			1 \\ \frac{1}{2c\sqrt{u_{33}}}
		\end{array}
		\right);
		W=
		\left(
		\begin{array}{c}
			W_{1} \\W_{2}
		\end{array}
		\right)
		=
		\left(
		\begin{array}{c}
			-\theta \\ 1
		\end{array}
		\right).
	\end{equation*}
	Furthermore, we get
	\begin{equation*}
		G_{n}(E_{33},n_{SN})=
		\left(
		\begin{array}{c}
			\frac{u(u-1)(u+m)}{(u+n)^{2}} \\ 0
		\end{array}
		\right)
		_{(E_{33},n_{SN})}
		=\left(
		\begin{array}{c}
			\frac{u_{33}(u_{33}-1)(u_{33}+m)}{(u_{33}+m_{2})^{2}} \\ 0
		\end{array}
		\right),
	\end{equation*}
	\begin{equation*}
		\begin{aligned}
			D^{2}G(E_{33},n_{SN})(V,V)&=
			\left(
			\begin{array}{c}
				\frac{{\partial }^{2}{G}_{1}}{\partial {u}^{2}}V_{1} ^{2}+\frac{{\partial }^{2}{G}_{1}}{\partial {u} \partial {v}}V_{1}V_{2}+\frac{{\partial }^{2}{G}_{1}}{\partial {v}^{2}}V_{2} ^{2}  \\
				\frac{{\partial }^{2}{G}_{2}}{\partial {u}^{2}}V_{1} ^{2}+\frac{{\partial }^{2}{G}_{2}}{\partial {u} \partial {v}}V_{1}V_{2}+\frac{{\partial }^{2}{G}_{2}}{\partial {v}^{2}}V_{2} ^{2}
			\end{array}
			\right)
			_{(E_{33},n_{SN})}                          \\
			&= \left(
			\begin{array}{c}
				\frac{-2u^{3}_{33}-6n_{SN}u^{2}_{33}-6n^{2}_{SN}u_{33}+(2m+2)n^{2}_{SN}+2mn_{SN}}{(u_{33}+n_{SN})^{3}}-\frac{1}{4cu_{33}} \\ -\frac{5\theta}{4cu_{33}}
			\end{array}
			\right)    ,                                    \\
		\end{aligned}
	\end{equation*}
	For $\kappa \neq0$, it's apparent that V and W satisfy
	\begin{equation*}
		W^{T}G_{n}(E_{33},n_{SN})=-	\frac{\theta u_{33}(u_{33}-1)(u_{33}+m)}{(u_{33}+m_{2})^{2}} \neq0,
	\end{equation*}
	\begin{equation*}
		W^{T}[D^{2}F(E_{21},c_{SN})(V,V)]=	\frac{\theta(2u^{3}_{33}+6n_{SN}u^{2}_{33}+6n^{2}_{SN}u_{33}-(2m+2)n^{2}_{SN}-2mn_{SN})}{(u_{33}+n_{SN})^{3}}-\frac{\theta}{cu_{33}}\neq0.
	\end{equation*}
	From Sotomayor's theorem in Perko \cite{Perko}, the saddle-node  bifurcation  appears at $E_{33}(u_{33},v_{33})$.
	
	\section*{Conflict of Interest}
The authors declare that they have no conflict of interest.

%\section*{Data Availability Statement}No data was used for research in this article.

\section*{Contributions}
 We declare that all the authors have same contributions to this paper.

\end{document}